\RequirePackage{fix-cm}
\documentclass[smallextended]{svjour3}       
\smartqed  
\usepackage{graphicx}
\usepackage{graphicx} 
\usepackage{amsmath} 
\usepackage{amssymb}
\usepackage{color}
\usepackage{subfigure}
\allowdisplaybreaks[4]
\begin{document}

\title{A Stabilized Physics Informed Neural Networks Method for Wave Equations}


\author{Yuling Jiao         \and
        Yuhui Liu            \and
        Jerry Zhijian Yang    \and
        Cheng Yuan
}


\institute{Yuling Jiao\at
School of Mathematics and Statistics, Wuhan University, Wuhan 430072, China\\
\email{yulingjiaomath@whu.edu.cn}
\and
Yuhui Liu\at
School of Mathematics and Statistics, Wuhan University, Wuhan 430072, China\\
\email{liu\_yuhui@whu.edu.cn}
\and
Jerry Zhijian Yang\at
School of Mathematics and Statistics, Wuhan University, Wuhan 430072, China\\
\email{zjyang.math@whu.edu.cn}
\and
Cheng Yuan\at
School of Mathematics and Statistics, Wuhan University, Wuhan 430072, China\\
\email{yuancheng@whu.edu.cn}
}

\date{Received: date / Accepted: date}

\maketitle

\begin{abstract}
  In this article, we propose a novel Stabilized Physics Informed Neural Networks method (SPINNs) for solving wave equations. In general, this method not only demonstrates theoretical convergence but also exhibits higher efficiency compared to the original PINNs. By replacing the $L^2$ norm with $H^1$ norm in the learning of initial condition and boundary condition, we theoretically proved that the error of solution can be upper bounded by the risk in SPINNs. Based on this, we decompose the error of SPINNs into approximation error, statistical error and optimization error. Furthermore, by applying the approximating theory of $ReLU^3$ networks and the learning theory on Rademacher complexity, covering number and pseudo-dimension of neural networks, we present a systematical non-asymptotic convergence analysis on our method, which shows that the error of SPINNs can be well controlled if the number of training samples, depth and width of the deep neural networks have been appropriately chosen. Two illustrative numerical examples on 1-dimensional and 2-dimensional wave equations demonstrate that SPINNs can achieve a faster and better convergence than classical PINNs method. 

\end{abstract}

\section{Introduction}
During the past few decades, numerical methods of Partial differential equations (PDEs) have been widely studied and applied in various fields of scientific computation \cite{brenner2007mathematical,ciarlet2002finite,Quarteroni2008Numerical,Thomas2013Numerical}. Among these, due to the central significance in solid mechanics, acoustics and electromagnetism, the numerical solution for wave equation attracts considerable attention, and a lot of work has been done to analyze the convergence rate, improve the solving efficiency and deal with practical problems such as boundary conditions. For many real problems with complex region, however, designing an efficient and accurate algorithms with practical absorbing boundary conditions is still difficult, especially for problems with irregular boundary. Furthermore, in high-dimensional case, many traditional methods may become even intractable due to the Curse of Dimensionality, which leads to an exponential increase in degree of freedom with the dimension of problem.

More recently, inspired by the great success of deep learning in fields of natural language processing and computational visions, solving PDEs with deep learning has become as a highly promising topic \cite{lagaris1998artificial,anitescu2019artificial,Berner2020Numerically,han2018solving,lu2021deepxde,sirignano2018dgm}. Several numerical schemes have been proposed to solve PDEs using neural networks, including the deep Ritz method (DRM) \cite{Weinan2017The}, physics-informed neural networks (PINNs) \cite{raissi2019physics}, weak adversarial neural networks(WANs) \cite{Yaohua2020weak} and their extensions \cite{jagtap2020conservative,npinns,fpinns}. Due to the simplicity and flexibility in its formulation, PINNs turns out to be the most concerned method. In the field of wave equations, researchers have successfully applied PINNs to the modeling of scattered acoustic fields \cite{wang2023acoustic}, including transcranial ultrasound wave \cite{wang2023physics} and seismic wave \cite{ding2023self}. In these works, all of the authors observed an interesting phenomenon that training PINNs without any boundary constraints may lead to a solution under absorbing boundary condition. In another word, the waves obtained by PINNs without boundary loss will be naturally absorbed at the boundary. This phenomenon, in fact, greatly improves the application value of PINNs in wave simulation, especially for inverse scattering problems. On the other hand, although PINNs have been widely used in the simulation of waves, a rigorous numerical analysis of PINNs for wave equations and more efficient training strategy are still needed.  

In this work, we propose the Stabilized Physics Informed Neural Networks (SPINNs) for simulation of waves. By replacing the $L^2$ norm in initial condition and boundary condition with $H^1$ norm, we obtain a stable PINNs method, in the sense that the error in solution can be upper bounded by the risk during training. It is worth mentioning that, in 2017 a similar idea called Sobolev Training has been proposed to improve the efficiency for regression \cite{czarnecki2017sobolev}. Later in \cite{son2021sobolev} and \cite{vlassis2021sobolev}, the authors generalized this idea to the training of PINNs, with applications to heat equation, Burgers' equation, Fokker-Planck equation and elasto-plasticity models. One main difference between our model and these works is that, we still use the $L^2$ norm, rather than $H^1$ norm for the residual and initial velocity in the loss of SPINNs. This designing, as we will demonstrate, turns out to be a sufficient condition to guarantee the stability, which also enables us to achieve a lower error with same and even less training samples. Furthermore, based on this stability, we firstly give a systematical convergence rate analysis of PINNs for wave equations. In general, our main contributions are summarized as follows:


\textbf{Main contributions}

$\bullet$ We propose a novel Stabilized Physics Informed Neural Networks method (SPINNs) for solving wave equations, for which we prove that the error in solution can be bounded by the training risk.

$\bullet$ We numerically show that SPINNs can achieve a faster and better convergence than original PINNs.

$\bullet$ We present a systematical convergence analysis of SPINNs. According to our result, once the network depth, width, as well as the number of training samples have been appropriately chosen, the error between the numerical solution $u_\phi$ from SPINNs and the exact solution $u^*$ can be arbitrarily small in the $H^1$ norm
\begin{align}
    \mathbb{E}_{\{X_n\}_{n=1}^N ,\{Y_m\}_{m=1}^M, \{T_k\}_{k=1}^K} \| \hat{u}_{\phi} - u^*\|_{H^1(\overline{\Omega_T})} \leq \varepsilon. \nonumber
\end{align}

The rest of this paper is organized as follows. In Section \ref{sec:method}, we describe the problem setting and introduce the SPINNs method. In Section \ref{sec:analysis}, we study the convergence rate of the SPINNs method for solving wave equations. In Section \ref{sec:experiment}, we present several numerical examples to illustrate the efficiency of SPINNs. Finally in Section \ref{sec:conclusion} the main conclusion will be provided.

\section{The Stabilized PINNs method}\label{sec:method}

In this section, we would introduce a stabilized PINNs (SPINNs) method for solving the wave equation. For completeness, we first list the following notations of neural networks and function spaces we will use. After that, the formulation of SPINNs will be presented.

\subsection{Preliminary}

Let $\mathcal{D} \in \mathbb{N}$, we would call function $f$ a neural network if it is implemented by:
\begin{align*}
f_0(x)& =  x, \\ 
f_l(x) & =  \rho_l(A_l f_{l-1}+b_l) ,\quad for \quad l =1,\cdots ,\mathcal{D}-1,  \\ 
f: & =  f_\mathcal{D}(x) =A_\mathcal{D}f_{\mathcal{D}-1} +b_\mathcal{D},  
\end{align*}
where $A_l = (a_{ij}^{(l)}) \in \mathbb{R}^{n_l \times n_{l-1}} $ ,
$b_l = (b_i^{(l)}) \in \mathbb{R}^{n_l}$.And $\rho_l: \mathbb{R}^{n_l} \rightarrow \mathbb{R}^{n_l}$, is the active function. The hyper-parameters $\mathcal{D}$ and $\mathcal{W}:= max\{N_l,l = 0,\cdots,\mathcal{D}\}$ are called the depth and the width of the network, respectively.
Let $\Phi$ be a set of activation functions and $X$ be a Banach space, the normed neural network function class can be defined as
\begin{align*}
    \mathcal{N}(\mathcal{D},\mathcal{W},\{\| \cdot \|_X, \mathcal{B} \},\Phi) := \{f: f \ \text{is implemented by a neural network with }  \\ 
    \text{depth } \mathcal{D} \text{ and width } \mathcal{W},\| f\|_X \leq \mathcal{B}, \rho_i^{(l)}\in \Phi\ \text{for each $i$ and $l$}\}.
\end{align*}
Next, we introduce several standard function spaces, including the continuous function space and Sobolev space:
\begin{align*}
    & C(\Omega) :=  \{\text{all the continuous functions defined on } \Omega \},\\ 
    & C^s(\Omega)  :=  \{ f: \Omega \rightarrow \mathbb{R}\ |\ D^{\alpha}f \in C(\Omega)\}, |\alpha| \leq s,\\ 
    & C(\overline{\Omega})  :=  \{\text{all the continuous functions defined on } \overline{\Omega} \}, \left\| f \right\|_{C(\overline{\Omega})} := \max_{x \in \overline{\Omega}}|f(x)|, \\
    & C^s(\overline{\Omega})  :=  \{ f: \Omega \rightarrow \mathbb{R}\ |\ D^{\alpha}f \in C(\overline{\Omega})\}, |\alpha| \leq s,\left\| f \right\|_{C^s(\overline{\Omega})} := \max_{x \in \overline{\Omega} , |\alpha| \leq s}|D^{\alpha}f(x)|, \\
    & L^p(\Omega)  :=  \left\{ f: \Omega \rightarrow \mathbb{R} | \int_\Omega \left|f\right|^p dx < \infty \right\},\\
    &\left\| f \right\|_{L^p(\Omega)} := \bigg[\int_\Omega |f|^p dx \bigg]^{1/p},\forall p \in [1,\infty),  \\
    & L^\infty(\Omega)  :=  \left\{ f: \Omega \rightarrow \mathbb{R}\ |\ \exists C >0\ s.t.|f|\leq C \ a.e. \right\},\\
    &\left\| f \right\|_{L^\infty(\Omega)} := \inf\{C\ |\ |f| \leq C \ a.e.\},  \\
    & H^s(\Omega)  :=   \left\{ f: \Omega \rightarrow \mathbb{R}\ |\  D^{\alpha}f \in L^2(\Omega),|\alpha| \leq s \right\},\\
    &\left\| f \right\|_{H^s(\Omega)} := \left(\sum_{|\alpha| \leq s} \left\| D^\alpha f \right\|_{L^2(\Omega)}^2\right)^{1/2}. 
\end{align*}

\subsection{Stabilized PINNs for wave equations}

Considering the following wave equation:

\begin{equation}\label{eq1}
\begin{cases}
u_{tt} - \Delta u = f , \quad\  \quad & (x,t) \in \Omega_T, \\
u(x,0) = \varphi (x) ,& x \in \Omega, \\
u_t (x,0) = \psi(x) ,& x \in \Omega, \\
u(x,t) = g(x,t) ,& x \in \partial \Omega, t \in [0,T],\\
\end{cases}
\end{equation}
where $\Delta u = \sum_{i=1}^d u_{x_i x_i}$, $\Omega = (0,1)^d$ and $\Omega_T = \Omega \times \left[0,T \right]$. 
Through this article, we would assume this problems defines a unique solution:
\begin{assumption}\label{strong solution}
Assume (\ref{eq1}) has a unique strong solution $u^* \in C^2(\Omega_T)\cap C(\overline{\Omega_T})$.
\end{assumption}

Without loss of generality, we further assume that $f,\varphi,\psi,g$ and their derivatives are $L^\infty$ bounded by a constant $\kappa$. Denote $\mathcal{B} \triangleq \max\{2 \| u^* \|_{C^2(\overline{\Omega_T})},\kappa^4 \}$.
Instead of solving problem \eqref{eq1} by traditional numerical methods, we turn to formulate \eqref{eq1} as a minimization problem on $C^2(\Omega_T)\cap C(\overline{\Omega_T})$, with the loss functional $\mathcal{L}$ being defined as:
\begin{align}\label{EPINN}
   \mathcal{L}(u) & :=  \| u_{tt} (x,t) - \Delta u(x,t) -f \|^2_{L^{2}(\Omega_T)} 
   + \| u(x,0) -  \varphi(x) \|^2_{H^{1}(\Omega)}  \nonumber\\ 
    &  +\| u_t (x,0) - \psi(x)\|^2_{L^{2}(\Omega)} 
    + \| u(x,t) -  g(x,t) \|^2_{H^{1}(\partial \Omega_T)}.
\end{align}

\begin{remark}
Different from the original loss in PINNs, we adopt $H^1$-norm in stead of $L^2$-norm in the learning of initial position and boundary condition. This modification, as we will demonstrate, offers advantages in both theoretical analysis and numerical computation.
\end{remark}

With assumption \ref{strong solution}, we know that $u^*$ is also the unique minimizer of loss functional $\mathcal{L}$ such that $\mathcal{L}(u^*)=0$. Let $|\Omega|$ and $|\partial\Omega|$ be the measure of $\Omega$ and $\partial\Omega$, namely, $|\Omega| := \int_\Omega 1 dx$, $|\partial\Omega| := \int_{\partial\Omega} 1 dx$ and $|T| := \int_0^T 1 dt$, then $\mathcal{L}(u)$ can be equivalently written as 
\begin{align}
    \mathcal{L}(u) 
    & = |\Omega||T| \mathbb{E}_{X \in U(\Omega),T \in U([ 0,T ])} \bigg( u_{tt}(X,T) -\Delta u(X,T) -f(X,T) \bigg) ^2\nonumber\\
    &+ |\Omega| \mathbb{E}_{X \in U(\Omega)} \bigg[\big(u(X,0) - \varphi(X)\big)^2 + \sum_{i=1}^d \big(u_{x_i}(X,0)- \varphi_{x_i} (X)\big)^2 \nonumber\\
    &+ \big(u_t(X,0) - \psi(X)\big)^2  \bigg] \nonumber\\
    & +|\partial \Omega||T| \mathbb{E}_{Y \in U(\partial\Omega),T \in U([ 0,T ])}\bigg[\big(u(Y,T)-g(Y,T)\big)^2 \nonumber\\
    &  +\big(u_t(Y,T)- g_t(Y,T)\big)^2 + \sum_{i=1}^d \big(u_{x_i}(Y,T) -  g_{x_i}(Y,T)\big)^2 \bigg] 
\end{align}
where $U(\Omega)$, $U(\partial\Omega)$, $U([0,T])$ are uniform distribution on $\Omega$, $\partial\Omega$ and $[0,T]$, respectively. To solve minimization of $\mathcal{L}(u)$ approximately, a Monte Carlo discrete version of $\mathcal{L}$ will be used:
\begin{align}
   \hat{\mathcal{L}}(u)
    & = \frac{|\Omega||T|}{N K} \sum_{n =1 }^{N} \sum_{k=1}^{K} \bigg( u_{tt}(X_n,T_k) -\Delta u(X_n,T_k) -f(X_n,T_k) \bigg) ^2\nonumber\\
    &+ \frac{|\Omega|}{N} \sum_{n =1 }^{N} \bigg[\big(u(X_n,0) - \varphi(X_n)\big)^2 + \sum_{i=1}^d\big(u_{x_i}(X_n,0)- \varphi_{x_i} (X_n)\big)^2 \nonumber\\
    &+ \big(u_t(X_n,0) - \psi(X_n)\big)^2  \bigg] \nonumber\\
    & +\frac{|\partial \Omega||T|}{M K} \sum_{m =1 }^{M} \sum_{k=1}^{K}\bigg[\big(u(Y_m,T_k)-g(Y_m,T_k)\big)^2 \nonumber\\
    &  +\big(u_t(Y_m,T_k)- g_t(Y_m,T_k)\big)^2 + \sum_{i=1}^d\big(u_{x_i}(Y_m,T_k) -  g_{x_i}(Y_m,T_k)\big)^2 \bigg] 
\end{align}
where $\{X_n\}_{n=1}^N$,$\{Y_m\}_{m=1}^M$ and $\{T_k\}_{k=1}^K$ are independent and identically distributed random samples according to the uniform distribution $U(\Omega)$, $U(\partial\Omega)$ and $U([0,T])$, respectively. With this approximation, we would solve the original problem \eqref{eq1} by using the empirical risk minimization:
\begin{align}\label{minimization}
    \hat{u}_{\phi} = \arg\min_{u_\phi \in \mathcal{P}} \hat{\mathcal{L}}(u_\phi),
\end{align}
where the admissible set $\mathcal{P}$ refers to the deep neural network function class parameterized by $\phi$. In this work, we will choose $\mathcal{P}$ as the $ReLU^3$ network function space, to ensure $\mathcal{P} \subset C^2(\Omega_T) \cap C(\overline{\Omega_T})$. More precisely,
$$
    \mathcal{P} = \mathcal{N}(\mathcal{D},\mathcal{W},\{ \| \cdot \|_{C^2(\overline{\Omega_T})},\mathcal{B}\} ,\{ReLU^3\}),\nonumber
$$
$\{\mathcal{D},\mathcal{W}\}$ will be given later to ensure the desired accuracy. The $ReLU^3$ activation function is defined by
\begin{align*}
    \rho(x) = 
    \begin{cases}
        x^3 , & x \geq 0,\\
        0, & others.
    \end{cases}
\end{align*}
In practical, the minimizer of problem \eqref{minimization} is usually obtained through some optimization algorithm $\mathcal{A}$. We would denote the minimizer as $u_{\phi_{\mathcal{A}}}$.

\section{Convergence analysis of SPINNs}\label{sec:analysis}

In this section, we will present a systematical error analysis of SPINNs for wave equations. To begin with, we first review some basic notations and theorem in the PDEs theory on wave equations. 

For wave equation, its total energy consists of two parts: kinetic energy $U$ and potential energy $V$, both of which can be expressed by multiple integrals,
$$ U = \frac{1}{2} \int_\Omega u_{t}^{2}dx,$$
$$ V = \frac{1}{2}\int_\Omega   \left(\sum_{i=1}^{d}u_{x_i}^2 \right)dx,$$
and their sum is called energy integral, the total energy of the wave equation (\ref{eq1}) excluding a constant factor is denoted as,
\begin{align}\label{eq2}
    E(t) = \int_\Omega \left(u_{t}^{2} + \sum_{i=1}^{d}u_{x_i}^2 \right)dx .
\end{align}

\begin{theorem}[Energy stability]
\label{Energy Inequalities}
    We denote $E_0(t) := \int_\Omega u^2 dx$, which stands for the square norm estimation of $u$. We have the energy inequality as below.
    
    \begin{align*}
        E(t) + E_0(t) & \leq C(T) (E(0) +  E_0(0)+\int_0^T \int_\Omega f^2dx dt  \\
        & + 2\int_0^T \int_{\partial\Omega} |u_t|\cdot \| \nabla u \| dsdt ) .
    \end{align*}
\end{theorem}

\begin{proof}
    See Appendix \ref{Energy_app} for details.
\end{proof}

\subsection{Risk decomposition}
By the definition of $u^*$ and $u_{\phi_{\mathcal{A}}}$, we can decompose the risk in SPINNs as
\begin{align*}
    \mathcal{L}(u_{\phi_\mathcal{A}}) - \mathcal{L} (u^*) = & \mathcal{L}(u_{\phi_\mathcal{A}}) - \hat{\mathcal{L}}(u_{\phi_\mathcal{A}}) + \hat{\mathcal{L}}(u_{\phi_\mathcal{A}}) - \hat{\mathcal{L}} (\hat{u}_{\phi})+ \hat{\mathcal{L}} (\hat{u}_{\phi})   \\
 & - \hat{\mathcal{L}}(\overline{u}) + \hat{\mathcal{L}}(\overline{u}) - \mathcal{L} (\overline{u}) + \mathcal{L} (\overline{u}) - \mathcal{L} (u^*)  \\	    
 = &\bigg[ \hat{\mathcal{L}}(\hat{u}_{\phi}) - \hat{\mathcal{L}}(\overline{u}) \bigg] + \bigg[ \mathcal{L}(u_{\phi_\mathcal{A}}) -  \hat{\mathcal{L}}(u_{\phi_\mathcal{A}})\bigg] + \bigg[ \hat{\mathcal{L}}(\overline{u}) -  \mathcal{L} (\overline{u})\bigg]   \\
 &+ \bigg[ \mathcal{L} (\overline{u})  - \mathcal{L} (u^*) \bigg] + \bigg[  \hat{\mathcal{L}}(u_{\phi_\mathcal{A}})- \hat{\mathcal{L}} (\hat{u}_{\phi})\bigg],  
\end{align*}
where $\overline{u}$ is an arbitrarily element in $\mathcal{P}$. Since $\hat{\mathcal{L}}(\hat{u}_\phi) - \hat{\mathcal{L}}(\overline{u}) \leq 0 $, and $u$ is an arbitrarily element in $\mathcal{P}$, we have:

\begin{align*}
    \mathcal{L}(u_{\phi_\mathcal{A}}) - \mathcal{L} (u^*) 
    \leq  
	\underbrace{2 \mathop{sup}_{u \in \mathcal{P}} \bigg| \mathcal{L}(u) -  \hat{\mathcal{L}}(u)\bigg|}_{\varepsilon_{sta}} 
	+\underbrace{\mathop{inf}_{u \in \mathcal{P}}\bigg|  \mathcal{L} (u)  - \mathcal{L} (u^*) \bigg|}_{\varepsilon_{app}}  
	+ \underbrace{\bigg[  \hat{\mathcal{L}}(u_{\phi_A})- \hat{\mathcal{L}} (\hat{u_{\phi}})\bigg] }_{\varepsilon_{opt}} 
\end{align*}
Thus, we have decomposed the total risk into approximation error ($\varepsilon_{app}$), statistical error ($\varepsilon_{sta}$) and optimization error ($\varepsilon_{opt}$). While the approximation error describes the expressive power of $ReLU^3$ network, the statistical error is caused by the discretization of the Monte Carlo method and the optimization error represents performance of the solver $\mathcal{A}$ we use. In this work,  we compromisely assume that the neural network can be well trained such that $\varepsilon_{opt} = 0$, and leave the optimization error as future study. In this case, it can be found that $\hat{u}_{\phi} = u_{\phi_{\mathcal{A}}}$.

\subsection{Lower bound of risk}
Next, based on the energy stability of wave equations, we shall present a lower bound of risk $\mathcal{L}(\hat{u}_{\phi}) -\mathcal{L}(u^*)$ in SPINNs. As we will demonstrate later, the risk can be arbitrary small if the network and sample complexity have been well chosen, and thus we can assume $\mathcal{L}(\hat{u}_{\phi}) -\mathcal{L}(u^*)<1$. Let $v = \hat{u}_{\phi} -u^*$ be the error between numerical solution and exact solution, we have
\begin{equation}\label{eq14}
\begin{cases}
v_{tt} - \Delta v  = (\hat{u}_{\phi})_{tt} - \Delta \hat{u}_{\phi} - f \triangleq \tilde{f} , \quad\  \quad & x \in \Omega ,t \in [0,T], \\
v(x,0) = \hat{u}_{\phi}(x, 0)- \varphi (x) ,& x \in \Omega, \\
v_t (x,0) = (\hat{u}_{\phi})_t(x , 0) - \psi(x) ,& x \in \Omega, \\
v(x,t) = \hat{u}_{\phi}(x, 0 ) - g(x,t) ,& x \in \partial \Omega , t \in [0,T],\\
\end{cases}
\end{equation}
and $\|v\|_{C^2(\overline{\Omega_T})}\leq \frac{3}{2} \mathcal{B}$. By applying theorem (\ref{Energy Inequalities}) to equation (\ref{eq14}), we obtain
\begin{align*}
     &\| \hat{u}_{\phi} -u^* \|^2_{H_1}   \\
    =& \int_\Omega \left(v_{t}^{2} + \sum_{i=1}^{d}v_{x_i}^2 \right)dx + \int_\Omega v^2 dx \\
    \leq &  C(T)\left( E(0) +  E_0(0)+\int_0^T \int_\Omega \tilde{f}^2dx dt + 2\int_0^T \int_{\partial\Omega} |v_t| \cdot \| \nabla v \| dsdt\right)\\
    \leq & C(T) \left( \int_\Omega(v_t(x,0)^2 +\sum_{i=1}^d v_{x_i}(x,0)^2) dx  + \int_ \Omega v(x,0)^2 dx \right.\\
    & \left.+ \int_{0}^T \int_\Omega((\hat{u}_{\phi})_ {tt}-\Delta\hat{u}_{\phi} -f)^2 dx dt + 3\sqrt{
    d}\mathcal{B}\int_0^T \int_{\partial\Omega} |v_t| dsdt  \right)\\
    \leq & C(T) \left( \int_\Omega(v_t(x,0)^2 +\sum_{i=1}^d v_{x_i}(x,0)^2) dx  + \int_ \Omega v(x,0)^2 dx \right.\\
    & \left.+ \int_{0}^T \int_\Omega((\hat{u}_{\phi})_ {tt}-\Delta\hat{u}_{\phi} -f)^2 dx dt + 3\sqrt{
    d}\mathcal{B}|\partial\Omega_T| (\int_0^T \int_{\partial\Omega} v_t^2 dsdt)^{\frac{1}{2}}  \right)\\
    \leq & C(T) \left(\| (\hat{u}_{\phi})_{tt} - \Delta \hat{u}_{\phi} -f \|^2_{L^{2}(\Omega_T)} +\| (\hat{u}_{\phi})_{t} (x,0) - \psi(x)\|^2_{L^{2}(\Omega)} \right. \\ 
    & \left.+\| \hat{u}_{\phi}(x,0) -  \varphi(x) \|^2_{H^{1}(\Omega)} + 3\sqrt{
    d}\mathcal{B}|\partial\Omega_T|\| \hat{u}_{\phi}(x,t) -  g(x,t) \|_{H^{1}(\partial \Omega_T)} \right)\\
    \leq & C(T) \left(\mathcal{L}(\hat{u}_{\phi}) +3\sqrt{d}\mathcal{B}|\partial\Omega_T|\mathcal{L}(\hat{u}_{\phi})^{\frac{1}{2}}\right)\\
    \leq & C(T)(1+3\sqrt{d}\mathcal{B}|\partial\Omega_T|) \mathcal{L}(\hat{u}_{\phi}) ^{\frac{1}{2}} \quad (\mathcal{L}(\hat{u}_{\phi})<1)\\
    \leq & \gamma \bigg(\mathcal{L}(\hat{u}_{\phi})-\mathcal{L}(u^*)\bigg)^{\frac{1}{2}},\quad (\mathcal{L}(u^*)=0)
\end{align*}
where we define $\gamma(T,d,\mathcal{B},|\partial\Omega_T|)\triangleq C(T)(1+3\sqrt{d}\mathcal{B}|\partial\Omega_T|)$. Combined this lower bound with previous risk decomposition, we can arrive at:
\begin{equation}\label{err_bound}
  \| \hat{u}_{\phi} -u^* \|^4_{H_1} \leq \gamma^2 (\varepsilon_{app} + \varepsilon_{sta}).
\end{equation}

\subsection{Approximation error}

By applying the following lemma proved in our previous work, we can get the upper bound of approximation error:
\begin{lemma}\label{thm_app}
    $\forall \overline{u} \in C^3(\overline{\Omega_T})$ and $\varepsilon > 0$, there exist a $ReLU^3$ network $u_\phi$ with depth $[log_2 d]+2$ and width $C(d,\| \overline{u} \|_{C^3(\overline{\Omega_T})})(\frac{1}{\varepsilon})^{d+1}$ such that $$\|\overline{u} - u_\phi \| _{C^2(\overline{\Omega_T})} \leq \varepsilon. $$
\end{lemma}

\begin{proof}
    A special case of Corollary 4.2 in \cite{AAM-39-239}.
\end{proof}

\begin{theorem}\label{thm_app_err}
    Under the  Assumption \ref{strong solution} and the condition that $u^* \in C^3(\overline{\Omega_T})$, for any $\varepsilon >0$, if we choose the following neural network function class:
    \begin{eqnarray*}
      \mathcal{P} &=& \mathcal{N}([log_2 d]+2 ,C(d,|\Omega|,|\partial\Omega|,|T|,\| u^* \|_{C^3(\overline{\Omega_T})}) (\frac{1}{\varepsilon})^{d+1}, \nonumber\\
       && \{\| \cdot \|_{C^2(\overline{\Omega_T})}, 2\| u^* \|_{C^2(\overline{\Omega_T})}\}, \{ReLU^3 \}  ), 
    \end{eqnarray*}
then the approximation error $\varepsilon_{app} \leq C(d,|\Omega_T|,\partial \Omega_T|) \varepsilon^2$.
\end{theorem}

\begin{proof}
    See Appendix \ref{Apperror_app} for details.
\end{proof}

\subsection{Statistical error}
The following theorem demonstrates that with sufficiently large sample complexity, the statistical error can be well controlled:
\begin{theorem}\label{statisticalerror}
    Let $\mathcal{D},\mathcal{W} \in \mathbb{N} , \mathcal{B} \in \mathbb{R}^+.$ For any $\varepsilon \geq 0 $, if the number of samples satisfy:
    \begin{equation*}
    \begin{cases}
        N = C(d,|\Omega|,\mathcal{B})\mathcal{D}^4\mathcal{W}^2(\mathcal{D}+log(\mathcal{W}))(\frac{1}{\varepsilon})^{2+\delta},\\
        K = C(d,|T|,\mathcal{B})\mathcal{D}^2 f_K (\mathcal{D},\mathcal{W})(\frac{1}{\varepsilon})^{k_1},\\
        M = C(d,|\partial\Omega|,\mathcal{B})f_M (\mathcal{D},\mathcal{W})(\frac{1}{\varepsilon})^{k_2},\\
    \end{cases}
\end{equation*}
where $f_K(\mathcal{D},\mathcal{W}) \geq 1$, $f_M(\mathcal{D},\mathcal{W}) \geq 1$, $\delta$ is an arbitrarily small number such that
\begin{equation*}
    \begin{cases}
        k_1 +k_2 = 2+\delta,\\
        f_k(\mathcal{D},\mathcal{W}) \cdot f_M(\mathcal{D},\mathcal{W}) = \mathcal{D}^2\mathcal{W}^2(\mathcal{D}+log(\mathcal{W})),
    \end{cases}
\end{equation*}
then we have:
\begin{align*}
    \mathbb{E}_{\{X_n\}_{n=1}^N ,\{Y_m\}_{m=1}^M, \{T_k\}_{k=1}^K} \sup_{u\in \mathcal{P}} | \mathcal{L}(u) -\hat{\mathcal{L}}(u)| \leq \varepsilon 
\end{align*}
\end{theorem}
\begin{proof}
    See Appendix \ref{Staerror_app} for details.
\end{proof}

\subsection{Convergence rate of SPINNs}
With the preparation in last two sections on the bounds of approximation and statistical errors, we will give the main results in this section.

\begin{theorem}
    Under the  Assumption \ref{strong solution}  and the condition that $u^* \in C^3(\overline{\Omega_T})$. For any $\varepsilon >0$, if we choose the parameterized neural network class
    \begin{eqnarray*}
      \mathcal{P} &=& \mathcal{N}([log_2 d]+2 ,C(d,|\Omega|,|\partial\Omega|,|T|,\| u^* \|_{C^3(\overline{\Omega_T})}) (\frac{1}{\varepsilon^2})^{d+1}, \nonumber\\
       && \{\| \cdot \|_{C^2(\overline{\Omega_T})}, 2\| u^* \|_{C^2(\overline{\Omega_T})}\}, \{ReLU^3 \}) 
    \end{eqnarray*}
    and let the number of samples be
     \begin{equation*}
    \begin{cases}
        N = C(d,|\Omega|,\mathcal{B})(\frac{1}{\varepsilon^4})^{d+3+\delta},\\
        K = C(d,|T|,\mathcal{B})(\frac{1}{\varepsilon^4})^{k_1 + \tilde{k}_1},\\
        M = C(d,|\partial\Omega|,\mathcal{B})(\frac{1}{\varepsilon^4})^{k_2 + \tilde{k}_2},\\
    \end{cases}
\end{equation*}
where $\tilde{k}_1,\tilde{k}_2 \geq 0$, $\delta$ is arbitrarily small such that
\begin{equation*}
    \begin{cases}
        k_1 +k_2 = 2 + \delta,\\
        \tilde{k}_1 +\tilde{k}_2 = d + 1,
    \end{cases}
\end{equation*}
then we have:
\begin{align*}
    \mathbb{E}_{\{X_n\}_{n=1}^N ,\{Y_m\}_{m=1}^M, \{T_k\}_{k=1}^K} \| \hat{u}_{\phi} - u^*\|_{H^1(\overline{\Omega_T})} \leq \varepsilon \nonumber
\end{align*}
\end{theorem}

\begin{proof}
     By theorem \ref{thm_app_err}, if we set the neural network function class as:
    \begin{eqnarray}
      \mathcal{P} &=& \mathcal{N}([log_2 d]+2 ,C(d,|\Omega|,|\partial\Omega|,|T|,\| u^* \|_{C^3(\overline{\Omega_T})}) (\frac{1}{\varepsilon^2})^{d+1}, \nonumber\\
       && \{\| \cdot \|_{C^2(\overline{\Omega_T})}, 2\| u^* \|_{C^2(\overline{\Omega_T})}\}, \{ReLU^3 \}) 
    \end{eqnarray}
     the approximation error can be arbitrarily small:
     \begin{align}\label{err_app_bound}
         \varepsilon_{app} \leq \frac{\varepsilon^4}{2\gamma^2}
     \end{align}
     Without loss of generality we assume that $\varepsilon$ is small enough such that
     \begin{align}
         \|\hat{u}_{\phi}\|_{C^2(\overline{\Omega_T})} \leq \| u^* - \hat{u}_{\phi}\|_{C^2(\overline{\Omega_T})} +\| u^* \|_{C^2(\overline{\Omega_T})} \leq 2\| u^* \|_{C^2(\overline{\Omega_T})}.\nonumber
     \end{align}
     By theorem \ref{statisticalerror}, when the number of samples be:
     \begin{equation}
        \begin{cases}
        N = C(d,|\Omega|,\mathcal{B})(\frac{1}{\varepsilon^4})^{d+3+\delta},\\
        K = C(d,|T|,\mathcal{B})(\frac{1}{\varepsilon^4})^{k_1 + \tilde{k}_1},\\
        M = C(d,|\partial\Omega|,\mathcal{B})(\frac{1}{\varepsilon^4})^{k_2 + \tilde{k}_2},\\
        \end{cases}
    \end{equation}
    where $\delta$ is an arbitrarily positive number and
    \begin{equation}
        \begin{cases}
        k_1 +k_2 = 2+ \delta,\\
        \tilde{k}_1 +\tilde{k}_2 = d + 1,
        \end{cases}
    \end{equation}
    we have:
    \begin{align}\label{err_sta_bound}
         \mathbb{E}_{\{X_n\}_{n=1}^N ,\{Y_m\}_{m=1}^M, \{T_k\}_{k=1}^K}\varepsilon_{sta} \leq \frac{\varepsilon^4}{2\gamma^2}
     \end{align}
     Combining \eqref{err_bound}, \eqref{err_app_bound} and \eqref{err_sta_bound} together, we get the final result:
     \begin{align}
         & \mathbb{E}_{\{X_n\}_{n=1}^N ,\{Y_m\}_{m=1}^M, \{T_k\}_{k=1}^K}\| \hat{u}_{\phi} - u^*\|_{H^1(\overline{\Omega_T})} \nonumber\\
         \leq &\mathbb{E}_{\{X_n\}_{n=1}^N ,\{Y_m\}_{m=1}^M, \{T_k\}_{k=1}^K}[\gamma^2\left(\mathcal{L}(\hat{u}_{\phi}) - \mathcal{L}(u^*) \right)]^{1/4} \nonumber \\
         \leq &\varepsilon .
     \end{align}
\end{proof}

\section{Numerical Experiments}\label{sec:experiment}

In this section, we would use SPINNs to solve wave equation in both one dimension and two dimension. 

\subsection{1D example}
Consider the following 1D wave equation on $\Omega = [-2, 2]$ from $t=0$ to $t= T=8$:
\begin{equation}\label{p1}
\begin{cases}
  u_{tt} = u_{xx}, \quad &x \in \Omega, \quad 0\leq t \leq 8,\\
  u(0,x) = 0, \quad  &x \in \Omega, \\
  u_t(0,x) = 0, \quad & x \in \Omega, \\
  u(t,-2) = sin(0.8\pi t),\quad & 0\leq t \leq 8, \\
  u(t,2) = 0,\quad & 0\leq t \leq 8.
\end{cases}
\end{equation}
For the training with SPINNs, we use a four-layer $ReLU^3$ network with 64 neurons in each layer to approximate the solution. We choose the Adam algorithm to implement the minimization, and the initial learning rate is set as 1E-3. 

\begin{figure}
\centering  
\subfigure{
\includegraphics[scale=0.35]{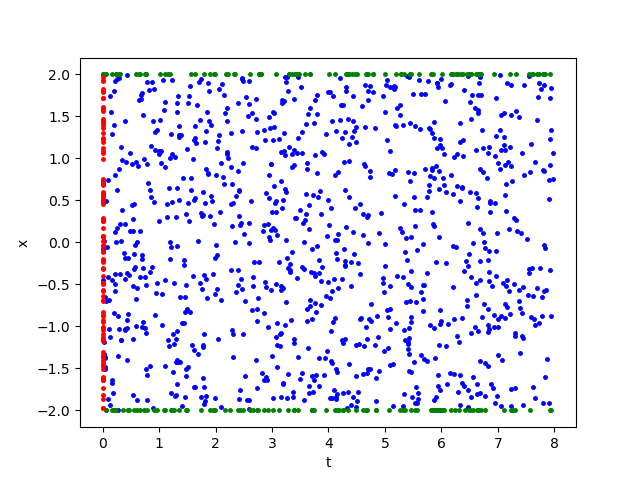}
}
\subfigure{
\includegraphics[scale=0.35]{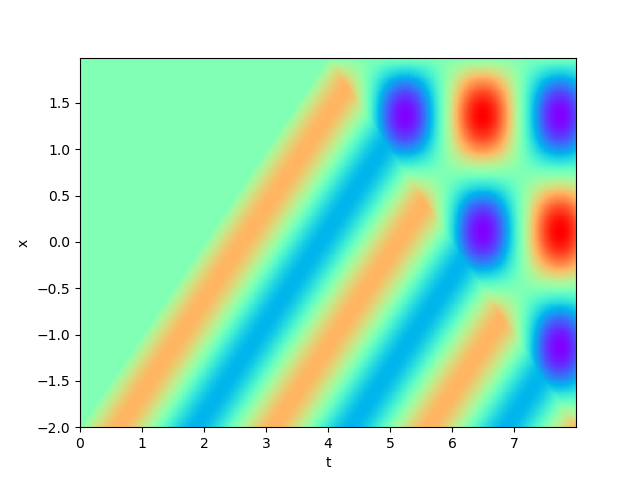}
}
\caption{Collection points (left) and reference solution (right) in $[0,T] \times \Omega$. The red, green and blue points stand for the initial points, boundary points and interior points correspondingly.}
\label{fig:exact_solution}     
\end{figure}

As for sample complexity, we train SPINNs with 10000 interior points, 500 boundary points (250 for each end) and 250 initial points in each epoch, all of which are sampled according to a uniform distribution (see (a) in Figure \ref{fig:exact_solution}). Further more, to obtain a better accuracy, we apply GAS method~\cite{jiao2023gas} as an adaptive sampling strategy. After every 250 epochs, we would adaptively add 600 inner points, 30 boundary points and 15 initial points based on a Gaussian mixture distribution. The GAS procedure will be repeated for 10 times. See \cite{jiao2023gas} for more details. For evaluation, we use the central difference method with a fine mesh ($dx = 0.01$, $dt = 0.009$) to obtain a reference solution $u_{cdm}$ (see (b) in Figure \ref{fig:exact_solution}) , with which we can calculate the following relative error by using numerical integration:
\begin{equation*}
  \text{Relative error} = \frac{\|u_{\phi} - u_{cdm}\|_{2}}{\|u_{cdm}\|_2}.
\end{equation*}

\begin{figure}
\centering  
\subfigure[$N_G = 2, u^{PINNs}$]{
\includegraphics[scale=0.33]{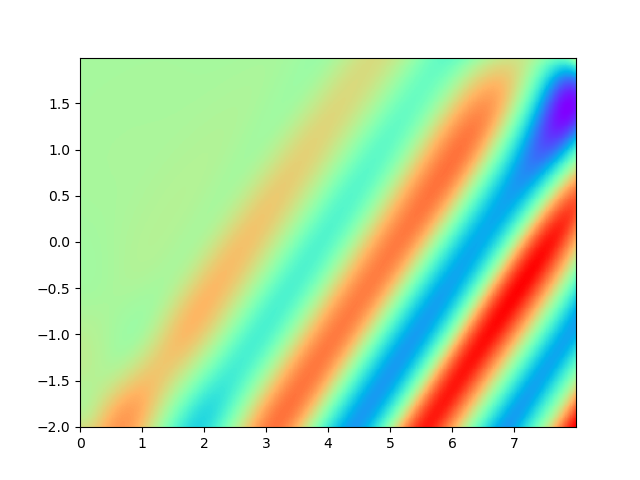}
}
\subfigure[$N_G = 2, u^{SPINNs}$]{
\includegraphics[scale=0.33]{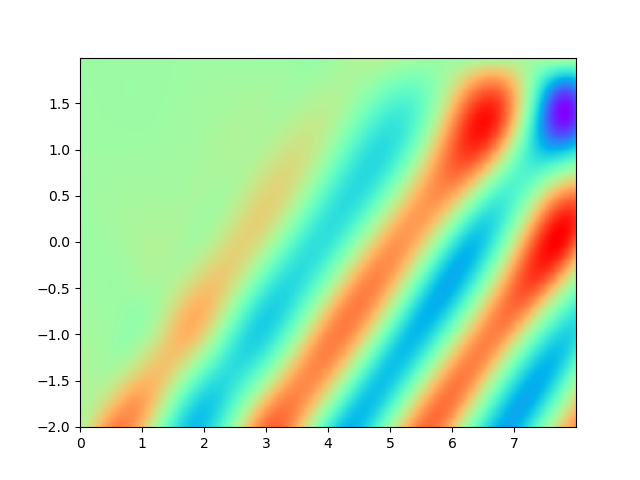}
}
\subfigure[$N_G = 5, u^{PINNs}$]{
\includegraphics[scale=0.33]{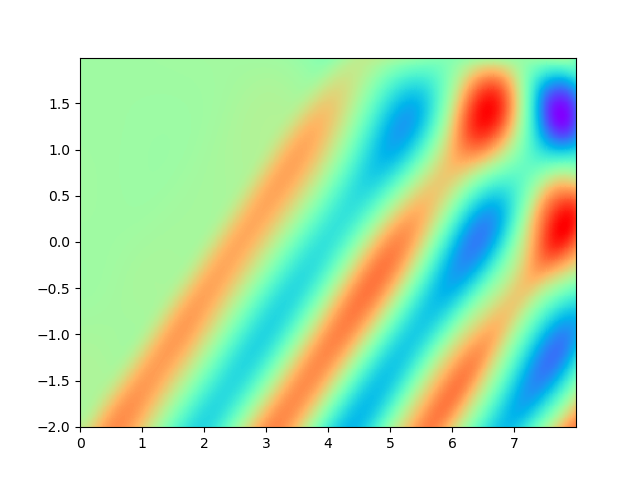}
}
\subfigure[$N_G = 5, u^{SPINNs}$]{
\includegraphics[scale=0.33]{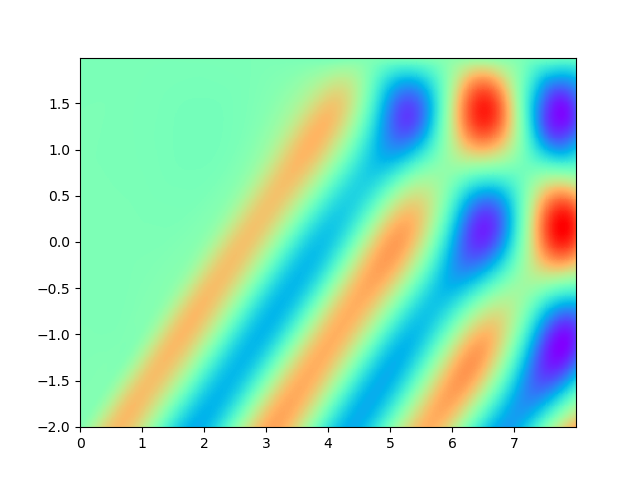}
}
\subfigure[$N_G = 8, u^{PINNs}$]{
\includegraphics[scale=0.33]{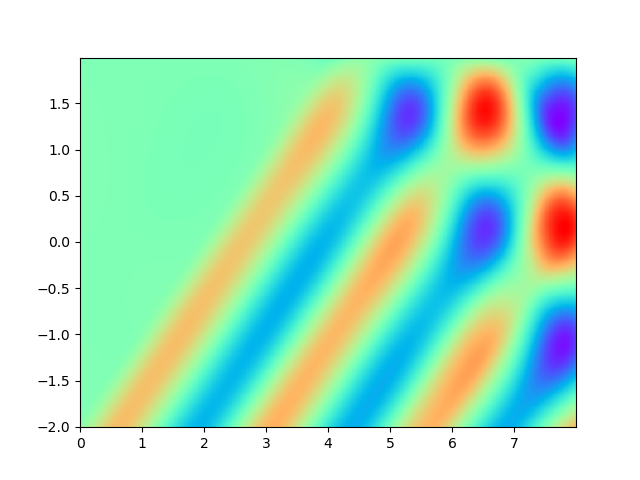}
}
\subfigure[$N_G = 8, u^{SPINNs}$]{
\includegraphics[scale=0.33]{fig/pinn_solution_8.png}
}
\subfigure[$N_G = 11, u^{PINNs}$]{
\includegraphics[scale=0.33]{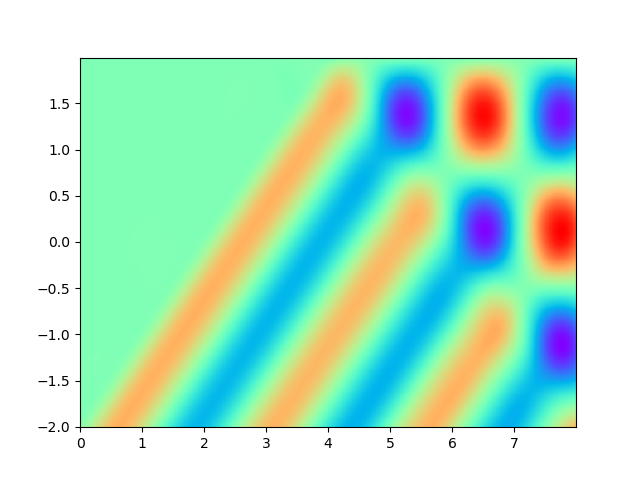}
}
\subfigure[$N_G = 11, u^{SPINNs}$]{
\includegraphics[scale=0.33]{fig/pinn_solution_11.png}
}
\caption{The numerical result of PINNs and SPINNs for equation \eqref{p1}.}    
\label{fig:1d_example}    
\end{figure}

\begin{figure}
  \centering
  \includegraphics[scale=0.5]{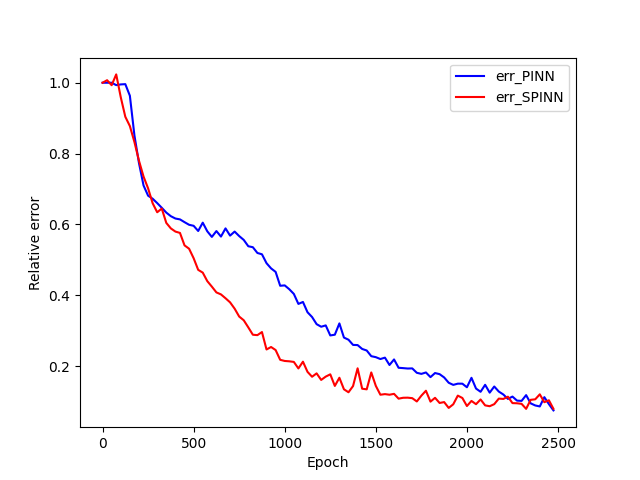}
  \caption{The relative error of PINNs and SPINNs.}\label{fig:rerr}
\end{figure}

Figure \ref{fig:1d_example} demonstrates the numerical result of PINNs and SPINNs for \eqref{p1} after $N_G$ times of adaptive sampling by GAS. As we can see, the SPINNs method converges faster than the classical PINNs, e.g., after 5 times of adaptive sampling, SPINNs have already captured all the six peaks of standing waves generated by the superposition of reflected wave and right-traveling waves. Furthermore, we present the relative error of PINNs and SPINNs in Figure \ref{fig:rerr}, which shows that our method can achieve a lower relative error at the early stage of training. On the other hand, with the times of adaptive sampling increasing, the classical PINNs can also arrive at a comparable accuracy, which has also been revealed by (g) in Figure \ref{fig:1d_example}. These results reflect the fact that training with SPINNs can speed up the convergence in solution, especially when the number of samples is relatively small.

\subsection{2D example}
Consider the following 2D wave equation on $\Omega = [-2, 2]^2$ from $t=0$ to $t= 1$:
\begin{equation}\label{p2}
\begin{cases}
  u_{tt} = u_{xx}+u_{yy}, \quad &x \in \Omega, \quad 0\leq t \leq 1,\\
  u(0,x,y) = sin(2\pi\sqrt{x^2+y^2}), \quad  &\sqrt{x^2+y^2} <1, \\
  u(0,x,y) = 0, \quad  &\sqrt{x^2+y^2} \geq 1, \\
  u_t(0,x,y) = 0, \quad & (x,y) \in \Omega, \\
  u(t,x,y) = 0,\quad & (x,y) \in \partial \Omega, 0\leq t \leq 1, \\
\end{cases}
\end{equation}
For the training with SPINNs, based on the experiment, we choose a three-layer $ReLU^3$ network with 512 neurons in each layer to approximate the solution. The optimization algorithm and the initial learning rate are kept as before.

\begin{figure}
\centering  
\subfigure[$t = 0, u_{\phi}$]{
\includegraphics[scale=0.33]{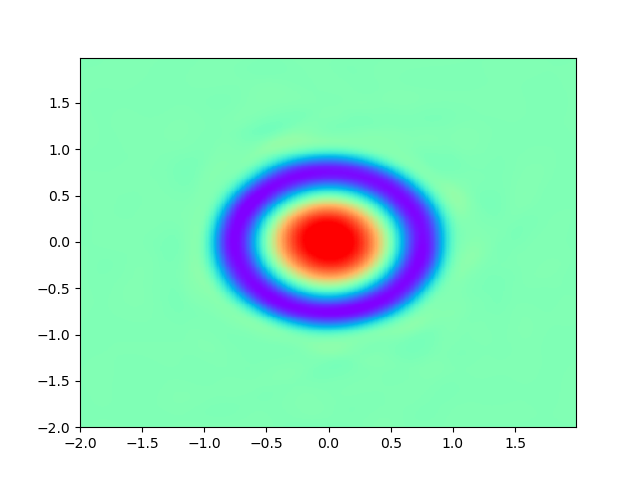}
}
\subfigure[$t = 0, |u_{cdm}-u_{\phi}|$]{
\includegraphics[scale=0.33]{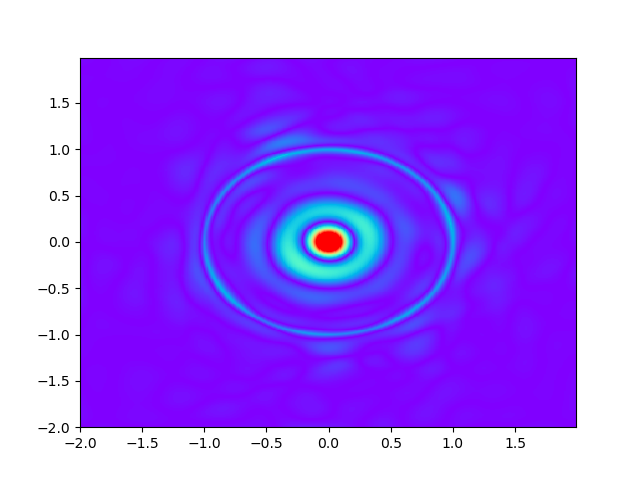}
}
\subfigure[$t = 0.1, u_{\phi}$]{
\includegraphics[scale=0.33]{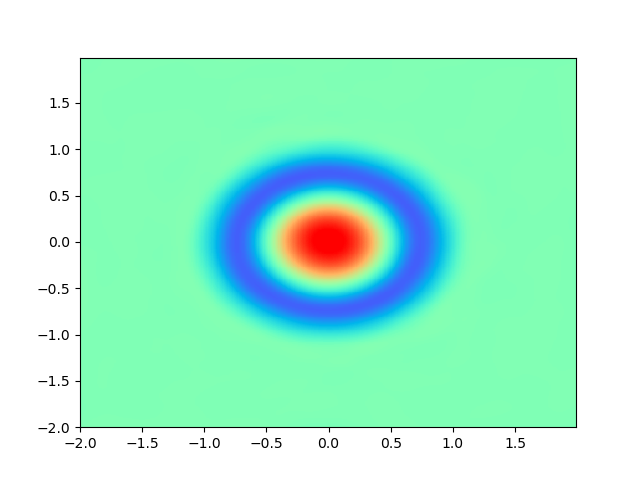}
}
\subfigure[$t = 0.1, |u_{cdm}-u_{\phi}|$]{
\includegraphics[scale=0.33]{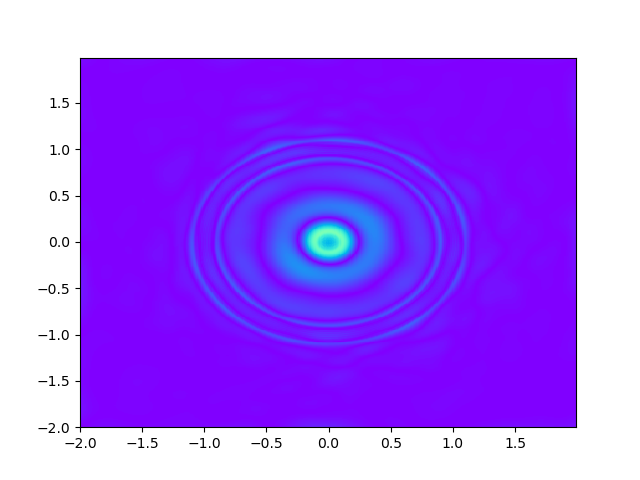}
}
\subfigure[$t = 0.2, u_{\phi}$]{
\includegraphics[scale=0.33]{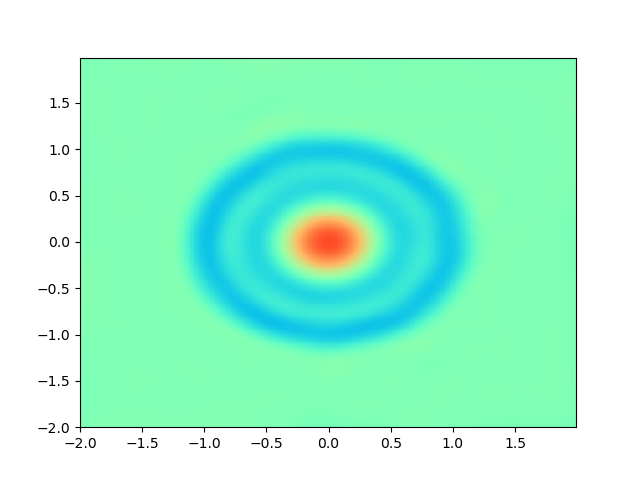}
}
\subfigure[$t = 0.2, |u_{cdm}-u_{\phi}|$]{
\includegraphics[scale=0.33]{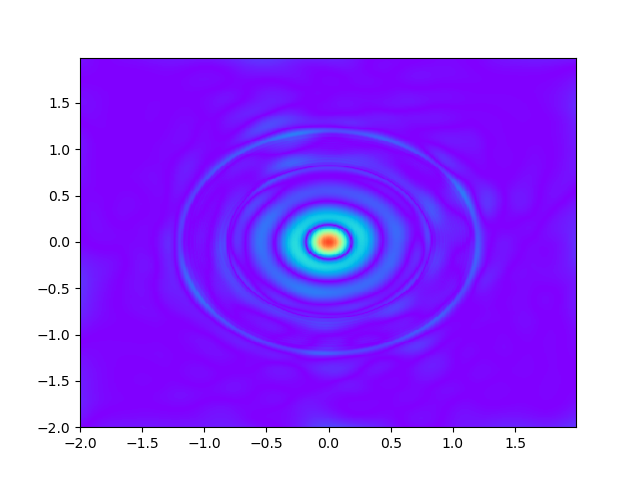}
}
\subfigure[$t = 0.3, u_{\phi}$]{
\includegraphics[scale=0.33]{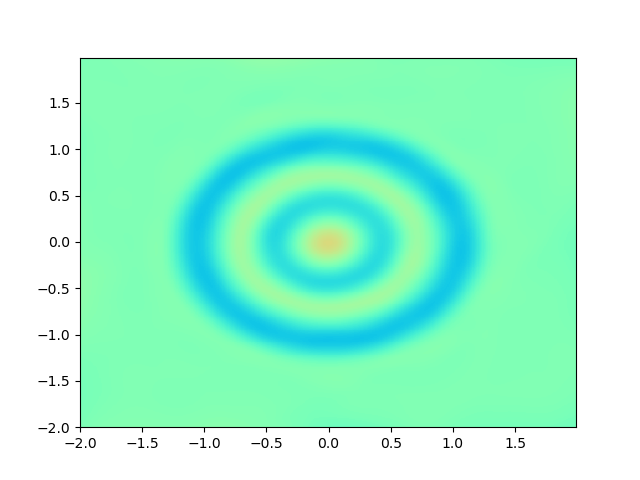}
}
\subfigure[$t = 0.3, |u_{cdm}-u_{\phi}|$]{
\includegraphics[scale=0.33]{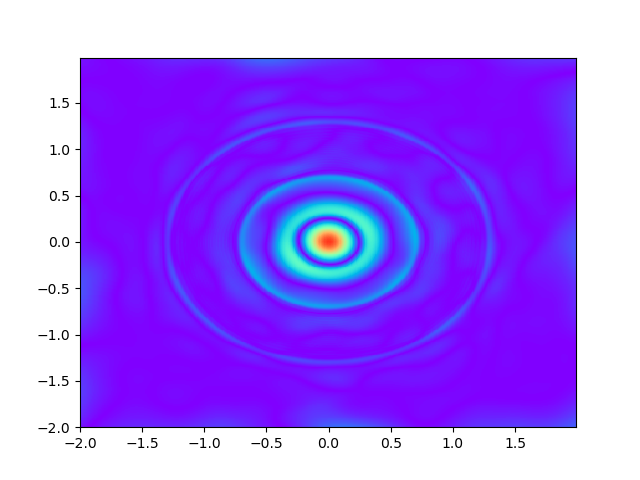}
}
\caption{The numerical result for \eqref{p2} after 10000 epochs training of PINNs.}    
\label{fig:2d_PINN}    
\end{figure}

\begin{figure}
\centering  
\subfigure[$t = 0, u_{\phi}$]{
\includegraphics[scale=0.33]{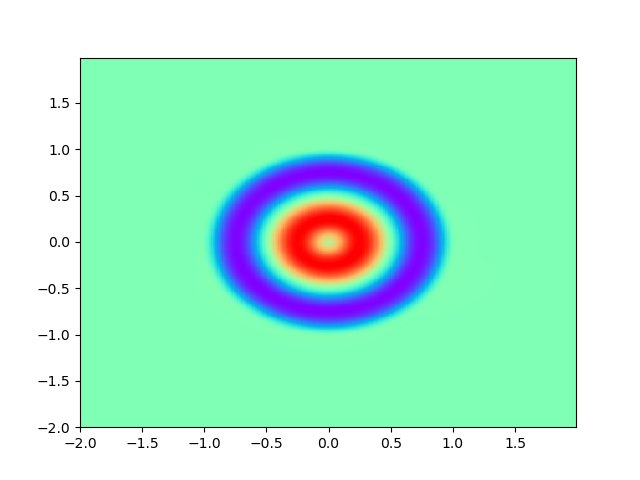}
}
\subfigure[$t = 0, |u_{cdm}-u_{\phi}|$]{
\includegraphics[scale=0.33]{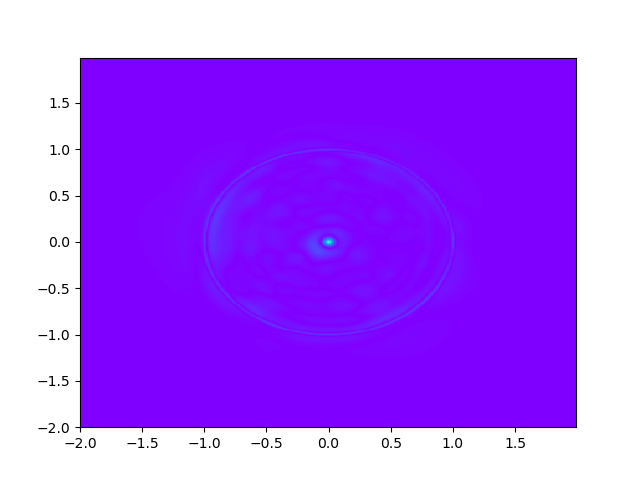}
}
\subfigure[$t = 0.1, u_{\phi}$]{
\includegraphics[scale=0.33]{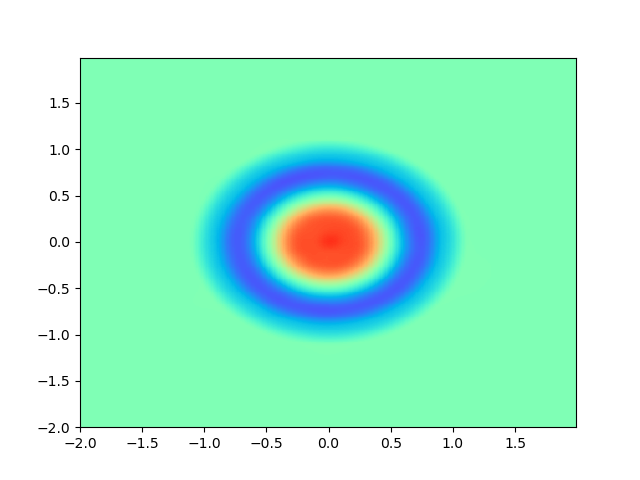}
}
\subfigure[$t = 0.1, |u_{cdm}-u_{\phi}|$]{
\includegraphics[scale=0.33]{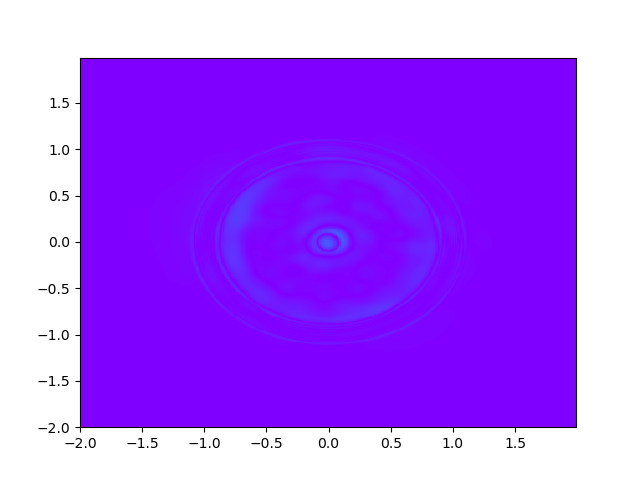}
}
\subfigure[$t = 0.2, u_{\phi}$]{
\includegraphics[scale=0.33]{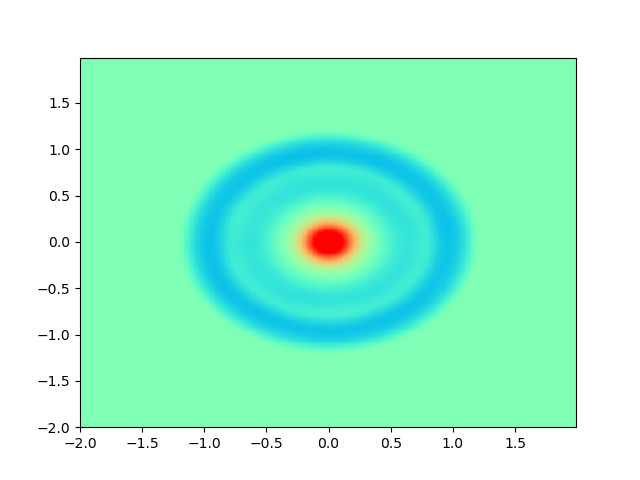}
}
\subfigure[$t = 0.2, |u_{cdm}-u_{\phi}|$]{
\includegraphics[scale=0.33]{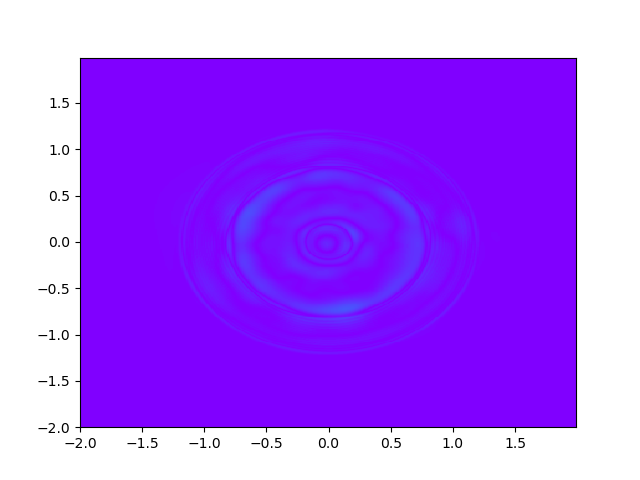}
}
\subfigure[$t = 0.3, u_{\phi}$]{
\includegraphics[scale=0.33]{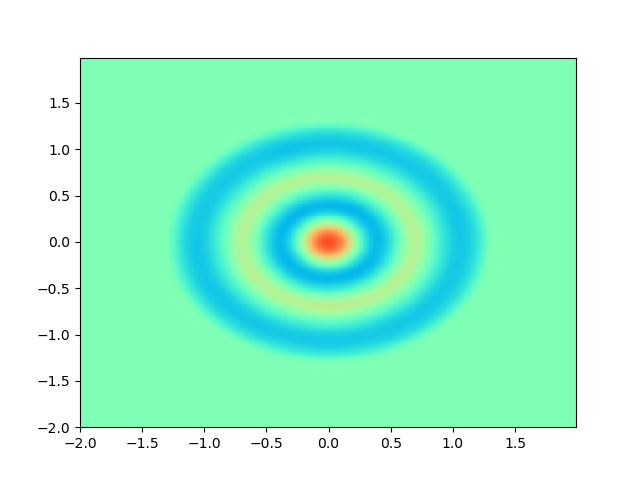}
}
\subfigure[$t = 0.3, |u_{cdm}-u_{\phi}|$]{
\includegraphics[scale=0.33]{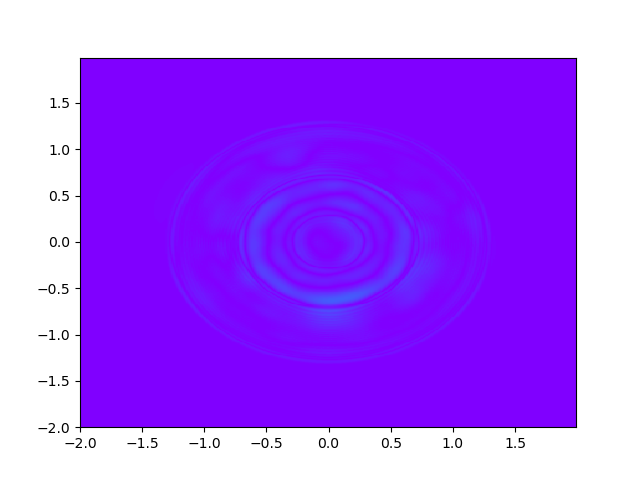}
}
\caption{The numerical result for \eqref{p2} after 10000 epochs training of SPINNs.}    
\label{fig:2d_EPINN}    
\end{figure}

As for sample complexity, we train SPINNs with 2000 interior points, 4000 boundary points (1000 for each edge) and 1000 initial points in each epoch, all of which are sampled according to a uniform distribution. For evaluation, we use the central difference method with a fine mesh ($dx=dy=0.01$, $dt = 0.004$) to obtain a reference solution $u_{cdm}$, with which we can calculate the pointwise absolute error $|u_{cdm}-u_{\phi}|$. As we can observe from Figure \ref{fig:2d_PINN} and Figure \ref{fig:2d_EPINN}, the SPINNs method can achieve a lower pointwise error, after training for same epochs. This superiority, in fact, can be understood as that learning with derivative information improves accuracy in the fitting of initial condition.

\section{Conclusion}\label{sec:conclusion}

In this work, we propose a stabilized physics informed neural networks method SPINNs for wave equations. With some numerical analysis, we rigorously prove SPINNs is a stable learning method, in which the solution error can be well controlled by the loss term. Based on this, a non-asymptotic convergence rate of SPINNs is presented, which provide people with a solid theoretical foundation to use it. Furthermore, by applying SPINNs to the simulation of two wave propagation problems, we numerically demonstrate that SPINNs can achieve a higher training accuracy and efficiency, especially when the number of samples is limited. On the other hand, how to extend this method to more difficult situations such as high dimensional problems and how to handel the optimization error in our convergence analysis are still needed to be studied. We will leave these topics as our future research.

\section*{Acknowledgement}
This work is supported by the National Key Research and Development Program
of China (No.2020YFA0714200), by the National Nature Science Foundation of China (No.12371441, No.12301558, No.12125103, No.12071362), and by the
Fundamental Research Funds for the Central Universities.

\section{Appendix}
\subsection{Appendix for energy integral of wave equations}\label{Energy_app}
According to the Gaussian formula, we have
\begin{align}\label{eq3}
    \int_\Omega \left(\sum_{i=1}^{d} (u_{t x_{i}}u_{x_i} )+ u_t (\sum_{i=1}^{d}u_{x_ix_i})\right)dx
    &= \int_\Omega (\nabla\cdot(u_t  \nabla u))dx \nonumber \\
    &= \int_{\partial \Omega} u_t \nabla u \cdot \mathbf{n} ds,
\end{align}
where $\mathbf{n}$ is the unit outer normal vector. Combine (\ref{eq2}) and (\ref{eq3}),we have
\begin{align}\label{eq4}
\frac{dE(t)}{dt}  & =  2\int_\Omega \left(u_t u_{tt}+ \sum_{i=1}^{d}u_{x_i} u_{x_{i}t}\right)dx \nonumber\\
                  & =  2\int_\Omega \left(u_t u_{tt} - u_t \sum_{i=1}^{d}u_{x_{i}x_{i}}\right)dx  + 2\int_{\partial\Omega} u_t \nabla u\cdot \mathbf{n}ds\nonumber\\
                  & = 2\int_\Omega u_t f dx +2\int_{\partial\Omega} u_t \nabla u\cdot \mathbf{n}ds\nonumber\\
                  & \leq \int_\Omega (u_t^2 + f^2 )dx +2\int_{\partial\Omega} |u_t|\cdot \| \nabla u \| ds.
\end{align}
Multiply both sides of the inequality (\ref{eq4}) by $e^{-t}$,
\begin{align}\label{eq5}
\frac{d(e^{-t}E(t))}{d t} \leq e^{-t} \left(\int_\Omega f^2 dx +2\int_{\partial\Omega} |u_t|\cdot \| \nabla u \| ds\right).
\end{align}
Then, integrating the equation (\ref{eq5}) from 0 to t,
\begin{align*}
E(t) \leq  e^t \left( E(0) + \int_0^t e^{- \tau} \int_\Omega f^2dx d\tau + 2\int_0^t e^{- \tau}\int_{\partial\Omega} |u_t|\cdot \| \nabla u \| dsd\tau \right)
\end{align*}
For any $t \in [0,T]$,
\begin{align}\label{eq7}
E(t) \leq C_1 \left(E(0) +  \int_0^T \int_\Omega f^2dx dt + 2\int_0^T \int_{\partial\Omega} |u_t|\cdot \| \nabla u \| dsdt\right),
\end{align}
$C_1$ is a constant that is only related to $T$.
Further we have 
\begin{align}\label{eq8}
\frac{dE_0(t)}{dt} = \int_\Omega 2uu_{t} dx \leq \int_\Omega u^2 dx +\int_\Omega u_t^2 dx \leq E_0(t) +E(t)
\end{align}
Multiply both sides of the above equation (\ref{eq8}) by $e^{-t}$,
\begin{align}\label{eq9}
\frac{d}{dt}(e^{-t}E_0(t)) \leq e^{-t}E(t)
\end{align}
Integrating the equation (\ref{eq9}) from 0 to t,
\begin{align}\label{eq10}
E_0(t) \leq e^{t}E_{0}(0)+e^{t} \int_0^t e^{-\tau}E(\tau) d \tau
\end{align}
For any $t \in [0,T]$,
\begin{align}\label{eq101}
E_0(t) \leq C_2\left(E_{0}(0)+E(t) \right),
\end{align}
$C_2$ is a constant that is only related to $T$.
Combine (\ref{eq7}) and (\ref{eq101}),we have,
\begin{align}\label{eq11}
E(t) +E_0(t) \leq C\left(E(0) +  E_0(0)+\int_0^T \int_\Omega f^2dx dt \right.\\ 
+ \left.2\int_0^T \int_{\partial\Omega} |u_t|\cdot \| \nabla u \| dsdt \right).
\end{align}
$C$ is a constant that is only related to $T$.

\subsection{Appendix for approximation error}\label{Apperror_app}

[Proof of theorem \ref{thm_app_err}] According to lemma \ref{thm_app}, we know that for $u^* \in C^3 (\overline{\Omega_T})$, $\varepsilon >0$, there exist a $ReLU^3$ network with depth $[log_2 d]+2$ and width $C(d,\| u^* \|_{C^3(\overline{\Omega_T})})(\frac{1}{\varepsilon})^{d+1}$, such that $\|v(x,t)\|_{C^2(\overline{\Omega_T})}=\| u^* -\hat{u}_{\phi} \|_{C^2(\overline{\Omega_T})} \leq  \varepsilon$. Hence,
\begin{align*}
    \varepsilon_{app} &\leq | \mathcal{L}(\hat{u}_{\phi})- \mathcal{L}(u^*) |  \\
                      & = \| (\hat{u}_{\phi})_{tt}(x,t) - \Delta \hat{u}_{\phi}(x,t) -f \|^2_{L^{2}(\Omega_T)} + \| \hat{u}_{\phi}(x,0) -  \varphi(x) \|^2_{H^{1}(\Omega)}  \\ 
                      &  + \| (\hat{u}_{\phi})_t (x,0) - \psi(x)\|^2_{L^{2}(\Omega)} 
                        + \| \hat{u}_{\phi}(x,t) -  g(x,t) \|^2_{H^{1}(\partial \Omega_T)}  \\
                      & \leq \| v_{tt} \|^2_{L^{2}(\Omega_T)} + \| \Delta v \|^2_{L^{2}(\Omega_T)} + \| v(x, 0)\|^2_{H^{1}(\Omega)} +\| v_t(x, 0)\|^2_{L^{2}(\Omega)} 
                        + \| v \|_{H^{1}(\partial \Omega_T)}^2 \\
                      & \leq C(d , |\Omega_T|, |\partial \Omega_T|) \cdot \varepsilon^2. 
\end{align*}

\subsection{Appendix for statistical error}\label{Staerror_app}

We will give the precise computation on the upper bounds of statistical error in this section. To begin with, we first introduce several basic concepts and results in learning theory.

\begin{definition}
    \textbf{(Rademacher complexity)} The Rademacher complexity of a set $A \subseteq \mathbb{R}^N$ is defined by
    \begin{align*}
        \Re(A) = \mathbb{E}_{\{\sigma_k\}_{k=1}^N} \Bigg[\ \sup_{a \in A} \frac{1}{N} \sum_{k=1}^{N} \sigma_k a_k \Bigg], 
    \end{align*}
    where $\{\sigma_k\}_{k=1}^N$ are $N$ i.i.d Rademacher variables with $p(\sigma_k = 1) = p(\sigma_k = -1) = \frac{1}{2}$. 

    Let $\Omega$ be a set and $\mathcal{F}$ be a function class which maps $\Omega$ to $\mathbb{R}$. Let $P$ be a probability
    distribution over $\Omega$ and $\{X_k\}_{k=1}^N$  be i.i.d. samples from $P$. The Rademacher complexity of $\mathcal{F}$ associated with distribution $P$ and sample size $N$ is defined by
    \begin{align*}
        \Re_{P,N}(\mathcal{F}) = \mathbb{E}_{\{X_k,\sigma_k\}_{k=1}^N} \Bigg[\ \sup_{u \in \mathcal{F}} \frac{1}{N} \sum_{k=1}^{N} \sigma_k u(X_k) \Bigg]. 
    \end{align*}
\end{definition}

\begin{lemma}
    Let $\Omega$ be a set and $P$ be a probability distribution over $\Omega$. Let $N \in \mathbb{N}$. Assume that $\omega: \Omega \rightarrow \mathbb{R}$ and $|\omega(x)| \leq \mathcal{B}$ for all $x \in \Omega$, then for any function class $\mathcal{F}$ mapping $\Omega$ to $\mathbb{R}$, there holds
    \begin{align*}
        \Re_{P,N}(\omega(x)\mathcal{F}) \leq 
        \mathcal{B}\Re_{P,N}(\mathcal{F}),
    \end{align*}
    where $\omega(x)\mathcal{F} := \{\overline{u}: \overline{u}(x) = \omega(x)u(x), u \in \mathcal{F}\}.$
\end{lemma}

\begin{proof}
    See \cite{jiao2022rate} for the proof.
\end{proof}

\begin{definition}
    \textbf{(Covering number)} Suppose that $W \subset \mathbb{R}$. For any $\varepsilon > 0$ ,let $V \subset \mathbb{R}^n$ be an $\varepsilon$-cover of $W$ with respect to the distance $d_\infty$, that is, for any $\omega \in W$, there exists a $v \in V$ such that $d_\infty(u,v) <\varepsilon$, where $d_\infty$ is defined by $d_\infty(u,v) := \max_{1\leq i \leq n} |u_i -v_i |$. The covering number $\mathcal{C}(\varepsilon,W,d_\infty)$ is defined to be the minimum cardinality among all $\varepsilon$-cover of $W$ with respect to the distance $d_\infty$.
\end{definition}

\begin{definition}
    \textbf{(Uniform covering number)} Suppose that $\mathcal{F}$ is a class of functions from $\Omega$ to $\mathbb{R}$. Given $n$ sample $Z_n =( Z_1,Z_2,\cdots Z_n) \in \Omega^n, \mathcal{F}|_{Z_n} \subset \mathbb{R}^n$ is defined by
    \begin{align*}
        \mathcal{F}|_{Z_n} \subset \mathbb{R}^n = 
        \{ (u(Z_1), u(Z_2), \cdots,u(Z_n)):u \in \mathcal{F} \}.
    \end{align*}
    The uniform covering number $\mathcal{C}_\infty(\varepsilon,\mathcal{F},n)$ is defined by
    \begin{align*}
        \mathcal{C}_\infty(\varepsilon,\mathcal{F},n) = \max_{Z_n \in \Omega^n}  \mathcal{C}(\varepsilon,\mathcal{F}|_{Z_n},d_\infty).
    \end{align*}
\end{definition}

\begin{lemma}\label{Rademacherbound}
    Let $\Omega$ be a set and $P$ be a probability distribution over $\Omega$. Let $N \in \mathbb{N}_{\geq 1}$,and $\mathcal{F}$ be a class of functions from $\Omega$ to $\mathbb{R}$ such that $0 \in \mathcal{F}$ and the diameter of $\mathcal{F}$ is less than $\mathcal{B}$,i.e.,$\| u\|_{L^\infty(\Omega)} \leq \mathcal{B}$,$\forall u \in \mathcal{F}$. Then
    \begin{align*}
        \Re_{P,N}(\mathcal{F}) \leq \inf_{0 < \delta <\mathcal{B}} \bigg(
        4 \delta +\frac{12}{\sqrt{N}}\int_{\delta}^\mathcal{B} \sqrt{log(2\mathcal{C}_\infty(\varepsilon,\mathcal{F},N))}d\varepsilon
        \bigg).
    \end{align*}
\end{lemma}

\begin{proof}
    This proof is base on the chaining method, see \cite{van1996weak}.
\end{proof}

\begin{definition}
    \textbf{(Pseudo-dimension)}
    Let $\mathcal{F}$ be a class of functions from $X$ to $\mathbb{R}$. Suppose that $S = \{x_1,x_2,\cdots,x_n\} \subset X $. We say that $S$ is pseudo-shattered by $\mathcal{F}$ if there exists $y_1,y_2,\cdots, y_n$ such that for any $b \in \{0,1\}^n$, there exists a $u \in \mathcal{F}$ satisfying
    \begin{align*}
        sign(u(x_i)-y_i) = b_i , i = 1,2,\cdots,n,
    \end{align*}
    and we say that $\{y_i\}_{i=1}^n$ witnesses the shattering. The pseudo-dimension of $\mathcal{F}$, denoted as $Pdim(\mathcal{F})$, is defined to be the maximum cardinality amoong all sets pseudo-shattered by $\mathcal{F}$.
\end{definition}
    
\begin{lemma}[\textbf{Theorem 12.2} in \cite{Anthony1999Neural}]\label{bound_cover}
    Let $\mathcal{F}$ be a class of real functions from a domain $X$ to the bounded interval $[0, \mathcal{B}]$.Let $\varepsilon >0$. Then 
    \begin{align*}
        \mathcal{C}_\infty(\varepsilon,\mathcal{F},n)\leq \sum_{i=1}^{Pdim(\mathcal{F})}\tbinom{n}{i} {\tbinom{\mathcal{B}}{\varepsilon}}^i ,
    \end{align*}
    which is less than $(\frac{en \mathcal{B}}{\varepsilon\cdot Pdim(\mathcal{F})})^{Pdim(\mathcal{F})}$ for $n \geq Pdim(\mathcal{F})$.
\end{lemma}

Next, to obtain the upper bound, we would decompose the statistical error into 24 terms by using triangle inequality:
\begin{align*}
    \mathbb{E}_{\{X_n\}_{n=1}^N ,\{Y_m\}_{m=1}^M, \{T_k\}_{k=1}^K} \sup_{u\in \mathcal{P}} | \mathcal{L}(u) -\hat{\mathcal{L}}(u)| \leq \\
    \sum_{j=1}^{24} \mathbb{E}_{\{X_n\}_{n=1}^N ,\{Y_m\}_{m=1}^M, \{T_k\}_{k=1}^K} \sup_{u\in \mathcal{P}} | \mathcal{L}_j(u) -\hat{\mathcal{L}}_j(u)|  
\end{align*}
where
\begin{align*}
    \mathcal{L}_1 & = |\Omega||T| \mathbb{E}_{X \sim U(\Omega) , T \sim U([0,T])} \left(u_{tt}(X,T)\right)^2 ,\\
    \mathcal{L}_2 & = |\Omega||T| \mathbb{E}_{X \sim U(\Omega) , T \sim U([0,T])} \left(\sum_{i=1}^d u_{x_i x_i}(X,T)\right)^2 ,\\
    \mathcal{L}_3 & = |\Omega||T| \mathbb{E}_{X \sim U(\Omega) , T \sim U([0,T])} \left(f(X,T)\right)^2 ,\\
    \mathcal{L}_4 & = -2|\Omega||T| \mathbb{E}_{X \sim U(\Omega) , T \sim U([0,T])} \left(\sum_{i=1}^d u_{tt}(X,T) u_{x_i x_i}(X,T)\right) ,\\
    \mathcal{L}_5 & = -2|\Omega||T| \mathbb{E}_{X \sim U(\Omega) , T \sim U([0,T])} \left(u_{tt}(X,T) f(X,T)\right) ,\\
    \mathcal{L}_6 & = 2|\Omega||T| \mathbb{E}_{X \sim U(\Omega) , T \sim U([0,T])} \left(\sum_{i=1}^d u_{x_i x_i}(X,T)f(X,T)\right) ,\\
    \mathcal{L}_7 & = |\Omega| \mathbb{E}_{X \sim U(\Omega)} \left(u(X,0)\right)^2 ,\\
    \mathcal{L}_8 & = |\Omega| \mathbb{E}_{X \sim U(\Omega)} \left(\varphi(X)\right)^2 ,\\
    \mathcal{L}_9 & = -2|\Omega| \mathbb{E}_{X \sim U(\Omega)} \left(u(X,0) \varphi(X)\right) ,\\
    \mathcal{L}_{10} & = |\Omega| \mathbb{E}_{X \sim U(\Omega)} \left(\sum_{i=1}^d u_{x_i}(X,0)^2\right) ,\\
    \mathcal{L}_{11} & = |\Omega| \mathbb{E}_{X \sim U(\Omega)} \left(\sum_{i=1}^d \varphi_{x_i}(X)^2\right) ,\\
    \mathcal{L}_{12} & = -2|\Omega| \mathbb{E}_{X \sim U(\Omega)} \left(\sum_{i=1}^d u_{x_i}(X,0)\varphi_{x_i}(X)\right) ,\\
    \mathcal{L}_{13} & = |\Omega| \mathbb{E}_{X \sim U(\Omega)} \left(u_t(X,0)\right)^2 ,\\
    \mathcal{L}_{14} & = |\Omega| \mathbb{E}_{X \sim U(\Omega)} \left(\psi(X)\right)^2 ,\\
    \mathcal{L}_{15} & = -2|\Omega| \mathbb{E}_{X \sim U(\Omega)} \left(u_t(X,0) \psi(X)\right) ,\\
    \mathcal{L}_{16} & = |\partial\Omega||T| \mathbb{E}_{Y \sim U(\partial\Omega) ,T \sim U([0,T])}\left( u(Y,T)\right)^2 ,\\
    \mathcal{L}_{17} & = |\partial\Omega||T| \mathbb{E}_{Y \sim U(\partial\Omega) ,T \sim U([0,T])}\left( g(Y,T)\right)^2 ,\\
    \mathcal{L}_{18} & = -2|\partial\Omega||T| \mathbb{E}_{Y \sim U(\partial\Omega) ,T \sim U([0,T])}\left( u(Y,T) g(Y,T)\right) ,\\
    \mathcal{L}_{19} & = |\partial\Omega||T| \mathbb{E}_{Y \sim U(\partial\Omega) ,T \sim U([0,T])}\left( u_t(Y,T)\right)^2 ,\\
    \mathcal{L}_{20} & = |\partial\Omega||T| \mathbb{E}_{Y \sim U(\partial\Omega) ,T \sim U([0,T])}\left( g_t(Y,T)\right)^2 ,\\
    \mathcal{L}_{21} & = -2|\partial\Omega||T| \mathbb{E}_{Y \sim U(\partial\Omega) ,T \sim U([0,T])}\left( u_t(Y,T) g_t(Y,T)\right) ,\\
    \mathcal{L}_{22} & = |\partial\Omega||T| \mathbb{E}_{Y \sim U(\partial\Omega) ,T \sim U([0,T])}\sum_{i=1}^d \left( u_{x_i}(Y,T)^2\right) ,\\
    \mathcal{L}_{23} & = |\partial\Omega||T| \mathbb{E}_{Y \sim U(\partial\Omega) ,T \sim U([0,T])}\sum_{i=1}^d \left( g_{x_i}(Y,T)^2\right) ,\\
    \mathcal{L}_{24} & = -2|\partial\Omega||T| \mathbb{E}_{Y \sim U(\partial\Omega) ,T \sim U([0,T])}\sum_{i=1}^d \left( u_{x_i}{x_i}(Y,T) g_{x_i}(X,T)\right),
\end{align*}
and $\hat{\mathcal{L}}_j(u)$ is the empirical version of $\mathcal{L}_j(u)$. The following lemma states that each of these 24 terms can be controlled by the corresponding Rademacher complexity.

\begin{lemma}\label{f-radermacher}
    Let ${\{X_n\}_{n=1}^N ,\{Y_m\}_{m=1}^M, \{T_k\}_{k=1}^K}$ be i.i.d samples from $U(\Omega),U(\partial \Omega), U([0,T])$, then we have
    \begin{align*}
        \mathbb{E}_{\{X_n\}_{n=1}^N ,\{Y_m\}_{m=1}^M, \{T_k\}_{k=1}^K}\sup_{u \in \mathcal{P} }\bigg| \mathcal{L}_j(u)- \hat{\mathcal{L}}_j(u) \bigg| \leq
        C(d,\mathcal{B})\Re_{\mathcal{U},N}(\mathcal{F}_j)  
    \end{align*}
    for $j = 1,2,\cdots ,24,$ where:
    \begin{align*}
        \mathcal{F}_{1} & = \{\pm f : \Omega_T \rightarrow \mathbb{R} | \quad \exists u\in \mathcal{P} \quad s.t. \quad f(x,t)= u_{tt}(x,t)^2 \},\\
        \mathcal{F}_{2} & = \{\pm f : \Omega_T \rightarrow \mathbb{R} | \quad \exists u\in \mathcal{P} \quad 1\leq i \leq j\leq d \quad s.t. \quad f(x,t)= u_{x_i x_i}(x,t)u_{x_j x_j}(x,t) \},\\
        \mathcal{F}_{4} & = \{\pm f : \Omega_T \rightarrow \mathbb{R} | \quad \exists u\in \mathcal{P} \quad 1\leq i \leq d \quad s.t. \quad f(x,t)= u_{tt}(x,t)u_{x_i x_i}(x,t) \},\\
        \mathcal{F}_{5} & = \{\pm f : \Omega_T \rightarrow \mathbb{R} | \quad \exists u\in \mathcal{P} \quad s.t. \quad f(x,t)= u_{tt}(x,t) \},\\
        \mathcal{F}_{6} & = \{\pm f : \Omega_T \rightarrow \mathbb{R} | \quad \exists u\in \mathcal{P} \quad 1\leq i \leq d \quad s.t. \quad f(x,t)= u_{x_i x_i}(x,t) \},\\
        \mathcal{F}_{7} & = \{\pm f : \Omega \rightarrow \mathbb{R} | \quad \exists u\in \mathcal{P} \quad s.t. \quad f(x)= u(x,0)^2 \},\\
        \mathcal{F}_{9} & = \{\pm f : \Omega \rightarrow \mathbb{R} | \quad \exists u\in \mathcal{P} \quad  s.t. \quad f(x)= u(x,0) \},\\
        \mathcal{F}_{10} & = \{\pm f : \Omega \rightarrow \mathbb{R} | \quad \exists u\in \mathcal{P} \quad 1\leq i \leq j\leq d \quad s.t. \quad f(x)= u_{x_i}(x,0)u_{x_j}(x,0) \},\\
        \mathcal{F}_{12} & = \{\pm f : \Omega \rightarrow \mathbb{R} | \quad \exists u\in \mathcal{P} \quad 1\leq i \leq d \quad s.t. \quad f(x,t)= u_{x_i}(x,0) \},\\
        \mathcal{F}_{13} & = \{\pm f : \Omega \rightarrow \mathbb{R} |\quad  \exists u\in \mathcal{P} \quad s.t. \quad f(x)= u_{t}(x,0)^2 \},\\
        \mathcal{F}_{15} & = \{\pm f : \Omega \rightarrow \mathbb{R} | \quad \exists u\in \mathcal{P} \quad s.t. \quad f(x)= u_{t}(x,0) \},\\
        \mathcal{F}_{16} & = \{\pm f : \partial \Omega_T \rightarrow \mathbb{R} | \quad \exists u\in \mathcal{P} \quad s.t. \quad f(x,t)= u(x,t)^2 |_{\partial \Omega} \},\\
        \mathcal{F}_{18} & = \{\pm f : \partial \Omega_T \rightarrow \mathbb{R} | \quad \exists u\in \mathcal{P} \quad s.t. \quad f(x,t)= u(x,t) |_{\partial \Omega} \},\\
        \mathcal{F}_{19} & = \{\pm f : \partial \Omega_T \rightarrow \mathbb{R} | \quad \exists u\in \mathcal{P} \quad s.t. \quad f(x,t)= u_t(x,t)^2 |_{\partial \Omega} \},\\
        \mathcal{F}_{21} & = \{\pm f : \partial \Omega_T \rightarrow \mathbb{R} | \quad \exists u\in \mathcal{P} \quad s.t. \quad f(x,t)= u_t(x,t) |_{\partial \Omega}  \},\\
        \mathcal{F}_{22} & = \{\pm f : \partial \Omega_T \rightarrow \mathbb{R} | \quad \exists u\in \mathcal{P} \quad 1\leq i \leq j \leq d\quad s.t. \quad f(x,t)= u_{x_i}(x,t)u_{x_j}(x,t) |_{\partial \Omega} \},\\
        \mathcal{F}_{24} & = \{\pm f : \partial \Omega_T \rightarrow \mathbb{R} | \quad \exists u\in \mathcal{P} \quad 1\leq i \leq d\quad s.t. \quad f(x,t)= u_{x_i}(x,t) |_{\partial \Omega} \}.
    \end{align*}
\end{lemma}

\begin{proof}
    The proof is based on the symmetrization technique, see \textbf{lemma 4.3} in \cite{jiao2022rate} for more details.
\end{proof}

\begin{lemma}\label{f-network}
   Let $\Phi=\{ReLU, ReLU^2, ReLU^3\}$. There holds
    \begin{align*}
        &\mathcal{F}_{1} \subset \mathcal{N}_{1} := \mathcal{N}(\mathcal{D}+5 ,(\mathcal{D}+2)(\mathcal{D}+4)\mathcal{W} ,\{\| \cdot \|_{C(\overline{\Omega_T})} , \mathcal{B}^2 \}, \Phi), \\
        &\mathcal{F}_{2} \subset \mathcal{N}_{2} := \mathcal{N}(\mathcal{D}+5 ,2(\mathcal{D}+2)(\mathcal{D}+4)\mathcal{W} ,\{\| \cdot \|_{C(\overline{\Omega_T})} , \mathcal{B}^2 \},  \Phi), \\
        &\mathcal{F}_{4} \subset \mathcal{N}_{4} := \mathcal{N}(\mathcal{D}+5 ,2(\mathcal{D}+2)(\mathcal{D}+4)\mathcal{W},\{\| \cdot \|_{C(\overline{\Omega_T})} , \mathcal{B}^2 \},  \Phi), \\
        &\mathcal{F}_{5} \subset \mathcal{N}_{5} := \mathcal{N}( \mathcal{D}+4, (\mathcal{D}+2)(\mathcal{D}+4)\mathcal{W},\{\| \cdot \|_{C(\overline{\Omega_T})} , \mathcal{B} \},  \Phi), \\
        &\mathcal{F}_{6} \subset \mathcal{N}_{6} := \mathcal{N}( \mathcal{D}+4, (\mathcal{D}+2)(\mathcal{D}+4)\mathcal{W},\{\| \cdot \|_{C(\overline{\Omega_T})} , \mathcal{B} \},  \Phi), \\
        &\mathcal{F}_{7} \subset \mathcal{N}_{7} := \mathcal{N}(\mathcal{D}+1) ,\mathcal{W} ,\{\| \cdot \|_{C(\overline{\Omega_T})} , \mathcal{B}^2 \},  \Phi), 
        \\
        &\mathcal{F}_{9} \subset \mathcal{N}_{9} := \mathcal{N}( \mathcal{D},\mathcal{W} ,\{\| \cdot \|_{C(\overline{\Omega_T})} , \mathcal{B} \},  \Phi), \\
        &\mathcal{F}_{10} \subset \mathcal{N}_{10} := \mathcal{N}( \mathcal{D}+2, 2(\mathcal{D}+2)\mathcal{W},\{\| \cdot \|_{C(\overline{\Omega_T})} , \mathcal{B}^2 \},  \Phi), \\
        &\mathcal{F}_{12} \subset \mathcal{N}_{12} := \mathcal{N}(\mathcal{D}+2, (\mathcal{D}+2)\mathcal{W},\{\| \cdot \|_{C(\overline{\Omega_T})} , \mathcal{B} \}, \Phi),\\
        &\mathcal{F}_{13} \subset \mathcal{N}_{13} := \mathcal{N}( \mathcal{D}+3, (\mathcal{D}+2)\mathcal{W},\{\| \cdot \|_{C(\overline{\Omega_T})} , \mathcal{B}^2 \},  \Phi),\\
        &\mathcal{F}_{15} \subset \mathcal{N}_{15} := \mathcal{N}( \mathcal{D}+2, (\mathcal{D}+2)\mathcal{W} ,\{\| \cdot \|_{C(\overline{\Omega_T})} , \mathcal{B} \},  \Phi),\\
        &\mathcal{F}_{16} \subset \mathcal{N}_{16} := \mathcal{N}( \mathcal{D}+1 ,\mathcal{W} ,\{\| \cdot \|_{C(\overline{\Omega_T})} , \mathcal{B}^2 \},  \Phi), \\
        &\mathcal{F}_{18} \subset \mathcal{N}_{18} := \mathcal{N}( \mathcal{D} ,\mathcal{W} ,\{\| \cdot \|_{C(\overline{\Omega_T})} , \mathcal{B} \}, \Phi), \\
        &\mathcal{F}_{19} \subset \mathcal{N}_{19} := \mathcal{N}( \mathcal{D}+3 ,(\mathcal{D}+2)\mathcal{W} ,\{\| \cdot \|_{C(\overline{\Omega_T})} , \mathcal{B}^2 \}, \Phi),\\
        &\mathcal{F}_{21} \subset \mathcal{N}_{21} := \mathcal{N}( \mathcal{D}+2 ,(\mathcal{D}+2)\mathcal{W}  ,\{\| \cdot \|_{C(\overline{\Omega_T})} , \mathcal{B} \},  \Phi), \\
        &\mathcal{F}_{22} \subset \mathcal{N}_{22} := \mathcal{N}( \mathcal{D}+3 ,2(\mathcal{D}+2)\mathcal{W} ,\{\| \cdot \|_{C(\overline{\Omega_T})} , \mathcal{B}^2 \}, \Phi),\\
        &\mathcal{F}_{24} \subset \mathcal{N}_{24} := \mathcal{N}( \mathcal{D}+2 ,(\mathcal{D}+2)\mathcal{W} ,\{\| \cdot \|_{C(\overline{\Omega_T})} , \mathcal{B} \}, \Phi). 
    \end{align*}
\end{lemma}

\begin{proof}
    The proof is an application of \textbf{proposition 4.2} in \cite{jiao2022rate}. Take $\mathcal{F}_{1}$ as an example, since $u\in \mathcal{P}$, we have $u_t\in\mathcal{N}(\mathcal{D}+2, (\mathcal{D}+2)\mathcal{W}, \{\|\cdot\|_{C^1(\overline{\Omega_T})}, \mathcal{B}\}, \Phi)$ and $u_{tt}\in \mathcal{N}(\mathcal{D}+4, (\mathcal{D}+2)(\mathcal{D}+4)\mathcal{W}, \{\|\cdot\|_{C(\overline{\Omega_T})}, \mathcal{B}\}, \Phi)$. Notice the square operation can be implemented as $x^2=ReLU^2(x) + ReLU^2(-x)$, thus we get that $u_{tt}^2 \in \mathcal{N}(\mathcal{D}+5, (\mathcal{D}+2)(\mathcal{D}+4)\mathcal{W}, \{\|\cdot\|_{C(\overline{\Omega_T})}, \mathcal{B}^2\}, \Phi)$. 
\end{proof}

\begin{lemma}[\textbf{Proposition 4.3} in \cite{jiao2022rate}]\label{pdim}
    For any $\mathcal{D},\mathcal{W} \in \mathbb{N}$,
    \begin{align*}
        Pdim(\mathcal{N}(\mathcal{D},\mathcal{W},\{ReLU,ReLU^2,ReLU^3\})) = \mathcal{O}(\mathcal{D}^2\mathcal{W}^2(\mathcal{D}+log\mathcal{W})).
    \end{align*}
\end{lemma}

Now we are ready to prove Theorem \ref{statisticalerror} on the statistical error.

\begin{proof}[The proof of \textbf{Theorem \ref{statisticalerror}}] According to lemma \ref{Rademacherbound} and lemma \ref{f-network}, 

$\bullet$ \quad For $i =1,2,4,5,6$, when the sample numbers $n=NK>Pdim(\mathcal{F}_i)$, we have

\begin{align}
    \Re_{P(\mathcal{X}),NK}(\mathcal{F}_i) &\leq \inf_{0<\delta <\mathcal{B}_i} \left( 4 \delta +\frac{12}{\sqrt{NK}} \int_{\delta}^{\mathcal{B}_i}\sqrt{log(2C_{\infty}(\varepsilon,\mathcal{F}_i,N))} d\varepsilon \right) \nonumber \\
    & \leq \inf_{0<\delta <\mathcal{B}_i} \left( 4 \delta +\frac{12}{\sqrt{NK}} \int_{\delta}^{\mathcal{B}_i}\sqrt{log \left( 2 \left(\frac{eNK\mathcal{B}_i}{\varepsilon\cdot Pdim(\mathcal{F}_i)}\right)^{Pdim(\mathcal{F}_i)}\right)}  d\varepsilon \right) \nonumber \\
    & \leq \inf_{0<\delta <\mathcal{B}_i} \left( 4 \delta + \frac{12\mathcal{B}_i}{\sqrt{NK}} +\frac{12}{\sqrt{NK}} \int_{\delta}^{\mathcal{B}_i}\sqrt{Pdim(\mathcal{F}_i)log \left( \frac{eNK\mathcal{B}_i}{\varepsilon\cdot Pdim(\mathcal{F}_i)}\right)}  d\varepsilon \right) \label{temp1}
\end{align}
Let $t = \sqrt{log\left( \frac{eNK\mathcal{B}_i}{\varepsilon\cdot Pim(\mathcal{F}_i)}\right)}$, then $\varepsilon = \frac{eNK\mathcal{B}_i}{Pdim(\mathcal{F}_i)} e^{-t^2}$. Denoting:
\begin{align*}
    t_1 = \sqrt{log\left( \frac{eNK\mathcal{B}_i}{\mathcal{B}_i\cdot Pdim(\mathcal{F}_i)} \right)} , t_2 = \sqrt{log\left( \frac{eNK\mathcal{B}_i}{\delta\cdot Pdim(\mathcal{F}_i)} \right)}
\end{align*}
we have:
\begin{align}
    &\int_{\delta}^{\mathcal{B}_i}\sqrt{log\left( \frac{eNK\mathcal{B}_i}{\varepsilon\cdot Pim(\mathcal{F}_i)}\right)}d\varepsilon \nonumber \\
    &= \frac{2eNK\mathcal{B}_i}{Pdim(\mathcal{F}_i)}\int_{t_1}^{t_2} t^2 e^{-t^2} dt \nonumber \\
    &= \frac{2eNK\mathcal{B}_i}{Pdim(\mathcal{F}_i)}\int_{t_1}^{t_2} t\left( \frac{-e^{-t^2}}{2}\right)^{'} dt \nonumber \\
    &= \frac{eNK\mathcal{B}_i}{Pdim(\mathcal{F}_i)} \left[ t_1 e^{-t_1^2} - t_2 e^{-t_2^2} +\int_{t_1}^{t_2} e^{-t^2} dt\right] \nonumber \\
    &\leq \frac{eNK\mathcal{B}_i}{Pdim(\mathcal{F}_i)} \left[ t_1 e^{-t_1^2} - t_2 e^{-t_2^2} + (t_2 - t_1) e^{-t_1^2}\right] \nonumber \\
    &\leq \frac{eNK\mathcal{B}_i}{Pdim(\mathcal{F}_i)} t_2 e^{-t_1^2} \nonumber \\
    &\leq \mathcal{B}_i \sqrt{log\left(\frac{eNK\mathcal{B}_i}{\delta\cdot Pdim(\mathcal{F}_i)}\right)} \label{temp2}
\end{align}    
Substitute \eqref{temp2} into \eqref {temp1} and choose $\delta = \mathcal{B}_i\left(\frac{Pdim(\mathcal{F}_i)}{NK}\right)^{\frac{1}{2}}\leq \mathcal{B}_i$, we have:
\begin{align}\label{bound NK}
    \Re_{P(\mathcal{X}),NK}(\mathcal{F}_i) &\leq \inf_{0<\delta <\mathcal{B}_i} \left( 4 \delta + \frac{12\mathcal{B}_i}{\sqrt{NK}} +\frac{12}{\sqrt{NK}} \int_{\delta}^{\mathcal{B}_i}\sqrt{Pdim(\mathcal{F}_i)log \left( \frac{eNK\mathcal{B}_i}{\varepsilon\cdot Pdim(\mathcal{F}_i)}\right)}  d\varepsilon \right)\nonumber\\
    & \leq \inf_{0<\delta <\mathcal{B}_i} \left( 4 \delta + \frac{12\mathcal{B}_i}{\sqrt{NK}} + \frac{12\mathcal{B}_i\sqrt{Pdim(\mathcal{F}_i)}}{\sqrt{NK}}\sqrt{log\left(\frac{eNK\mathcal{B}_i}{\delta \cdot Pdim(\mathcal{F}_i)}\right)}\right)\nonumber\\
    &\leq 28 \sqrt{\frac{3}{2}} \mathcal{B}_i \left( \frac{Pdim(\mathcal{F}_i)}{NK}\right)^{\frac{1}{2}}\sqrt{log\left(\frac{eNK}{Pdim(\mathcal{F}_i)} \right)}\nonumber\\
    &\leq 28 \sqrt{\frac{3}{2}} \mathcal{B}_i \left( \frac{Pdim(\mathcal{N}_i)}{NK}\right)^{\frac{1}{2}}\sqrt{log\left(\frac{eNK}{Pdim(\mathcal{N}_i)} \right)}\nonumber\\
    &\leq 28 \sqrt{\frac{3}{2}} max\{\mathcal{B},\mathcal{B}^2\} \left(\frac{\mathcal{H}_1}{NK}\right)^{\frac{1}{2}}\sqrt{log\left(\frac{eNK}{\mathcal{H}_1}\right)}
\end{align}
where $$\mathcal{H}_1 = C_1 (\mathcal{D}+2)^2(\mathcal{D}+4)^2(\mathcal{D}+5)^2\mathcal{W}^2[(\mathcal{D}+5)+log(2(\mathcal{D}+2)(\mathcal{D}+4)\mathcal{W})].$$ The last step above is due to lemma \ref{pdim}.

$\bullet$ \quad For $i= 7,9,10,12,13,15$, when the sample numbers $n=N>Pdim(\mathcal{F}_i)$, we can similarly prove that
\begin{align}\label{bound N}
    \Re_{\mathcal{U}(\Omega),N}(\mathcal{F}_i) \leq 28 \sqrt{\frac{3}{2}} max\{\mathcal{B},\mathcal{B}^2\} \left(\frac{\mathcal{H}_2}{N}\right)^{\frac{1}{2}}\sqrt{log\left(\frac{eN}{\mathcal{H}_2} \right)}
\end{align}
where $$\mathcal{H}_2 = C_2 (\mathcal{D}+2)^2(\mathcal{D}+3)^2\mathcal{W}^2[(\mathcal{D}+3)+log(2(\mathcal{D}+2)\mathcal{W})].$$

$\bullet$ \quad For $i= 16,18,19,21,22,24$, when the sample numbers $n=MK>Pdim(\mathcal{F}_i)$, we have
\begin{align}\label{bound MK}
    \Re_{\mathcal{U}(\partial \Omega_T),MK}(\mathcal{F}_i) \leq 28 \sqrt{\frac{3}{2}} max\{\mathcal{B},\mathcal{B}^2\} \left(\frac{\mathcal{H}_3}{MK}\right)^{\frac{1}{2}}\sqrt{log\left(\frac{eMK}{\mathcal{H}_3} \right)}
\end{align}
where $$\mathcal{H}_3 = C_3 (\mathcal{D}+2)^2(\mathcal{D}+3)^2\mathcal{W}^2[(\mathcal{D}+3)+log(2(\mathcal{D}+2)\mathcal{W})].$$

$\bullet$ \quad For $i =3$, we have,
\begin{align*}
    &\mathbb{E}_{\{X_n\}_{n=1}^N , \{T_k\}_{k=1}^K} \bigg| \mathcal{L}_3  - \hat{\mathcal{L}}_3\bigg|\\
    &  =  |\Omega| |T| \mathbb{E}_{\{X_n\}_{n=1}^N , \{T_k\}_{k=1}^K} \bigg| \mathbb{E}_{X \sim U(\Omega) , T \sim U([0,T])} f(X,T)^2 - \frac{1}{NK}  \sum_{n=1}^{N} \sum_{k =1 }^{K} f(X_n,T_k)^2 \bigg|\\
    &  \leq  |\Omega| |T| \sqrt{\mathbb{E}_{\{X_n\}_{n=1}^N , \{T_k\}_{k=1}^K}\bigg| \mathbb{E}_{X \sim U(\Omega) , T \sim U([0,T])} f(X,T)^2 - \frac{1}{NK}  \sum_{n=1}^{N} \sum_{k =1 }^{K} f(X_n,T_k)^2 \bigg|^2}  \\
    &  = \frac{ |\Omega| |T| }{NK} \sqrt{\mathbb{E}_{\{X_n\}_{n=1}^N , \{T_k\}_{k=1}^K} \sum_{n=1}^{N} \sum_{k =1 }^{K} \bigg| \mathbb{E}_{X \sim U(\Omega) , T \sim U([0,T])} f(X,T)^2 - f(X_n,T_k)^2 \bigg|^2} \\
    &  = \frac{ |\Omega| |T| }{NK} \sqrt{NK \cdot\mathbb{E}_{X_1 \sim U(\Omega) , T_1 \sim U([0,T])} \bigg| \mathbb{E}_{X \sim U(\Omega) , T \sim U([0,T])} f(X,T)^2 - f(X_1,T_1)^2 \bigg|^2} \\
    & = \frac{ |\Omega| |T| }{NK} \sqrt{NK}\sigma(f(X,T)^2)\\
    & = |\Omega| |T| \frac{\sigma(f(X,T)^2)}{\sqrt{NK}} 
\end{align*}
where $\sigma(f(X,T))$ is the standard deviation of $f(X,T)$. With the bound of $f$, we can further obtain,

\begin{align*}
    \sigma^2(f(X,T)^2) & = \mathbb{E}\left((f(X,T)^2)^2\right) - \left(\mathbb{E}(f(X,T)^2)\right)^2\\
    & \leq \frac{1}{|\Omega||T|}\left(\int_{\Omega_T} (f(X,T)^2)^2 dX dT\right)\\
    & \leq \frac{1}{|\Omega||T|}\left(\int_{\Omega_T} \kappa^4 dX dT\right)\\
    & \leq \frac{1}{|\Omega||T|} \cdot |\Omega||T| \mathcal{B}\\
    & = \mathcal{B}
\end{align*}
then we have,
\begin{align*}
    & \mathbb{E}_{X \sim U(\Omega) , T \sim U([0,T])} \bigg| \mathcal{L}_3  - \hat{\mathcal{L}}_3\bigg| \leq  |\Omega| |T| \sqrt{\frac{\mathcal{B}}{NK}} .
\end{align*}

$\bullet$  \quad Similarly, for $i= 8,14, 17, 20$,
\begin{align*}
    & \mathbb{E}_{X \sim U(\Omega) , T \sim U([0,T])} \bigg| \mathcal{L}_8  - \hat{\mathcal{L}}_8\bigg|  \leq |\Omega| \sqrt{\frac{\mathcal{B}}{N}} ,\\
    & \mathbb{E}_{X \sim U(\Omega)} \bigg| \mathcal{L}_{14}  - \hat{\mathcal{L}}_{14}\bigg|  \leq |\Omega| \sqrt{\frac{\mathcal{B}}{N}},\\
    & \mathbb{E}_{Y \sim U(\partial\Omega) , T \sim U([0,T])} \bigg| \mathcal{L}_{17}  - \hat{\mathcal{L}}_{17}\bigg|  \leq |\partial \Omega| |T| \sqrt{\frac{\mathcal{B}}{MK}}  ,\\
    & \mathbb{E}_{Y \sim U(\partial\Omega) , T \sim U([0,T])} \bigg| \mathcal{L}_{20}  - \hat{\mathcal{L}}_{20}\bigg|  \leq |\partial \Omega| |T| \sqrt{\frac{\mathcal{B}}{MK}} .
\end{align*}

$\bullet$ \quad For $i =11 $,
\begin{align*}
    \sigma^2(\sum_{i=1}^{d} \varphi_{x_i}(x)^2) &= \mathbb{E}\left((\sum_{i=1}^{d} \varphi_{x_i}(x)^2)^2\right) - \left(\mathbb{E}(\sum_{i=1}^{d} \varphi_{x_i}(x)^2)\right)^2\\
    & \leq \frac{1}{|\Omega|}\left[\int_{\Omega} \left(\sum_{i=1}^{d} \varphi_{x_i}(x)^2\right)^2 dx\right]  \\
    & \leq \frac{1}{|\Omega|}\left[\int_{\Omega}( d \kappa^2)^2 dx \right]\\
    & \leq d^2 \mathcal{B}
\end{align*}
then we have,
\begin{align*}
    & \mathbb{E}_{X \sim U(\Omega) } \bigg| \mathcal{L}_{11}  - \hat{\mathcal{L}}_{11}\bigg| \leq d|\Omega| \sqrt{\frac{\mathcal{B}}{N}}.
\end{align*}

$\bullet$  \quad Similarly, for $i= 23$,
\begin{align*}
    & \mathbb{E}_{Y \sim U(\partial\Omega) , T \sim U([0,T])} \bigg| \mathcal{L}_{23}  - \hat{\mathcal{L}}_{23}\bigg|  \leq d|\partial \Omega||T| \sqrt{\frac{\mathcal{B}}{MK}}.
\end{align*}

Hence, we have,
\begin{align*}
    &\mathbb{E}_{\{X_n\}_{n=1}^N ,\{Y_m\}_{m=1}^M, \{T_k\}_{k=1}^K} \sup_{u\in \mathcal{P}} | \mathcal{L}(u) -\hat{\mathcal{L}}(u)|  \\
    &\leq\sum_{j=1}^{24} \mathbb{E}_{\{X ,Y, T\}} \sup_{u\in \mathcal{P}} | \mathcal{L}_j(u) -\hat{\mathcal{L}}_j(u)|  \\
    &\leq 28 \sqrt{\frac{3}{2}} max\{\mathcal{B},\mathcal{B}^2\}\left( 5|\Omega| |T|C_1 \left(\frac{\mathcal{H}_1}{NK}\right)^{\frac{1}{2}}\sqrt{log\left(\frac{eNK}{\mathcal{H}_1}\right)} +  \right.\\
    &\left.6|\Omega| C_2 \left(\frac{\mathcal{H}_2}{N}\right)^{\frac{1}{2}}\sqrt{log\left(\frac{eN}{\mathcal{H}_2}\right)} + 6|\partial \Omega| |T|C_3 \left(\frac{\mathcal{H}_3}{MK}\right)^{\frac{1}{2}}\sqrt{log\left(\frac{eMK}{\mathcal{H}_3}\right)} \right)\\
    & \left( |\Omega| |T|\sqrt{\frac{\mathcal{B}}{NK}} + (2+d)|\Omega|\sqrt{\frac{\mathcal{B}}{N}} +(2+d)|\partial \Omega| |T|\sqrt{\frac{\mathcal{B}}{MK}}\right),
\end{align*}
where $C_1,C_2,C_3$ are constants associated with dimensionality $d$ and bound $\mathcal{B}$.
Hence, for any $\varepsilon \geq 0$, if the number of samples satisfies:
\begin{equation*}
    \begin{cases}
        N &= C(d,|\Omega|,\mathcal{B})\mathcal{D}^4\mathcal{W}^2(\mathcal{D}+log(\mathcal{W}))(\frac{1}{\varepsilon})^{2+\delta},\\
        K &= C(d,|T|,\mathcal{B})\mathcal{D}^2 f_K (\mathcal{D},\mathcal{W})(\frac{1}{\varepsilon})^{k_1},\\
        M &= C(d,|\partial\Omega|,\mathcal{B})f_M (\mathcal{D},\mathcal{W})(\frac{1}{\varepsilon})^{k_2},\\
    \end{cases}
\end{equation*}

where:
\begin{equation*}
    \begin{cases}
        k_1 +k_2 = 2+\delta,\\
        f_k(\mathcal{D},\mathcal{W}) \cdot f_M(\mathcal{D},\mathcal{W}) = \mathcal{D}^2\mathcal{W}^2(\mathcal{D}+log(\mathcal{W})).
    \end{cases}
\end{equation*}
with restriction $f_K(\mathcal{D},\mathcal{W}) \geq 1$, $f_M(\mathcal{D},\mathcal{W}) \geq 1$, and $\delta$ is arbitrarily small. Then we have:
\begin{align*}
    \mathbb{E}_{\{X_n\}_{n=1}^N ,\{Y_m\}_{m=1}^M, \{T_k\}_{k=1}^K} \sup_{u\in \mathcal{P}} | \mathcal{L}(u) -\hat{\mathcal{L}}(u)| \leq \varepsilon. 
\end{align*}

\end{proof}

    \bibliography{ref.bib}

\begin{thebibliography}{10}

\bibitem{brenner2007mathematical}
Susanne Brenner and Ridgway Scott.
\newblock {\em The mathematical theory of finite element methods}, volume~15.
\newblock Springer Science \& Business Media, 2007.

\bibitem{ciarlet2002finite}
Philippe~G Ciarlet.
\newblock {\em The finite element method for elliptic problems}.
\newblock SIAM, 2002.

\bibitem{Quarteroni2008Numerical}
A.~Quarteroni and A.~Valli.
\newblock {\em Numerical Approximation of Partial Differential Equations},
  volume~23.
\newblock Springer Science \& Business Media, 2008.

\bibitem{Thomas2013Numerical}
J.W. Thomas.
\newblock {\em Numerical Partial Differential Equations: Finite Difference
  Methods}, volume~22.
\newblock Springer Science \& Business Media, 2013.

\bibitem{lagaris1998artificial}
Isaac~E Lagaris, Aristidis Likas, and Dimitrios~I Fotiadis.
\newblock Artificial neural networks for solving ordinary and partial
  differential equations.
\newblock {\em IEEE transactions on neural networks}, 9(5):987--1000, 1998.

\bibitem{anitescu2019artificial}
Cosmin Anitescu, Elena Atroshchenko, Naif Alajlan, and Timon Rabczuk.
\newblock Artificial neural network methods for the solution of second order
  boundary value problems.
\newblock {\em Computers, Materials \& Continua}, 59(1):345--359, 2019.

\bibitem{Berner2020Numerically}
Julius Berner, Markus Dablander, and Philipp Grohs.
\newblock Numerically solving parametric families of high-dimensional
  kolmogorov partial differential equations via deep learning.
\newblock In {\em Advances in Neural Information Processing Systems},
  volume~33, pages 16615--16627. Curran Associates, Inc., 2020.

\bibitem{han2018solving}
Jiequn Han, Arnulf Jentzen, and E~Weinan.
\newblock Solving high-dimensional partial differential equations using deep
  learning.
\newblock {\em Proceedings of the National Academy of Sciences},
  115(34):8505--8510, 2018.

\bibitem{lu2021deepxde}
Lu~Lu, Xuhui Meng, Zhiping Mao, and George~Em Karniadakis.
\newblock Deepxde: A deep learning library for solving differential equations.
\newblock {\em SIAM Review}, 63(1):208--228, 2021.

\bibitem{sirignano2018dgm}
Justin Sirignano and Konstantinos Spiliopoulos.
\newblock Dgm: A deep learning algorithm for solving partial differential
  equations.
\newblock {\em Journal of computational physics}, 375:1339--1364, 2018.

\bibitem{Weinan2017The}
E.~Weinan and Ting Yu.
\newblock The deep ritz method: A deep learning-based numerical algorithm for
  solving variational problems.
\newblock {\em Communications in Mathematics and Statistics}, 6(1):1--12, 2017.

\bibitem{raissi2019physics}
Maziar Raissi, Paris Perdikaris, and George~E Karniadakis.
\newblock Physics-informed neural networks: A deep learning framework for
  solving forward and inverse problems involving nonlinear partial differential
  equations.
\newblock {\em Journal of Computational Physics}, 378:686--707, 2019.

\bibitem{Yaohua2020weak}
Yaohua Zang, Gang Bao, Xiaojing Ye, and Haomin Zhou.
\newblock Weak adversarial networks for high-dimensional partial differential
  equations.
\newblock {\em Journal of Computational Physics}, 411:109409, 2020.

\bibitem{jagtap2020conservative}
Ameya~D Jagtap, Ehsan Kharazmi, and George~Em Karniadakis.
\newblock Conservative physics-informed neural networks on discrete domains for
  conservation laws: Applications to forward and inverse problems.
\newblock {\em Computer Methods in Applied Mechanics and Engineering},
  365:113028, 2020.

\bibitem{npinns}
G.~Pang, M.~D'Elia, M.~Parks, and G.E. Karniadakis.
\newblock npinns: Nonlocal physics-informed neural networks for a parametrized
  nonlocal universal laplacian operator. algorithms and applications.
\newblock {\em Journal of Computational Physics}, 422:109760, 2020.

\bibitem{fpinns}
Guofei Pang, Lu~Lu, and George~Em Karniadakis.
\newblock fpinns: Fractional physics-informed neural networks.
\newblock {\em SIAM Journal on Scientific Computing}, 41(4):A2603--A2626, 2019.

\bibitem{wang2023acoustic}
Hao Wang, Jian Li, Linfeng Wang, Lin Liang, Zhoumo Zeng, and Yang Liu.
\newblock On acoustic fields of complex scatters based on physics-informed
  neural networks.
\newblock {\em Ultrasonics}, 128:106872, 2023.

\bibitem{wang2023physics}
Linfeng Wang, Hao Wang, Lin Liang, Jian Li, Zhoumo Zeng, and Yang Liu.
\newblock Physics-informed neural networks for transcranial ultrasound wave
  propagation.
\newblock {\em Ultrasonics}, 132:107026, 2023.

\bibitem{ding2023self}
Yi~Ding, Su~Chen, Xiaojun Li, Suyang Wang, Shaokai Luan, and Hao Sun.
\newblock Self-adaptive physics-driven deep learning for seismic wave modeling
  in complex topography.
\newblock {\em Engineering Applications of Artificial Intelligence},
  123:106425, 2023.

\bibitem{czarnecki2017sobolev}
Wojciech~M Czarnecki, Simon Osindero, Max Jaderberg, Grzegorz Swirszcz, and
  Razvan Pascanu.
\newblock Sobolev training for neural networks.
\newblock {\em Advances in neural information processing systems}, 30, 2017.

\bibitem{son2021sobolev}
Hwijae Son, Jin~Woo Jang, Woo~Jin Han, and Hyung~Ju Hwang.
\newblock Sobolev training for physics informed neural networks.
\newblock {\em Communications in Mathematical Sciences}, 21:1679--1705, 2013.

\bibitem{vlassis2021sobolev}
Nikolaos~N Vlassis and WaiChing Sun.
\newblock Sobolev training of thermodynamic-informed neural networks for
  interpretable elasto-plasticity models with level set hardening.
\newblock {\em Computer Methods in Applied Mechanics and Engineering},
  377:113695, 2021.

\bibitem{AAM-39-239}
Yuling Jiao, Xiliang Lu, Jerry~Zhijian Yang, Cheng Yuan, and Pingwen Zhang.
\newblock Improved analysis of pinns: Alleviate the cod for compositional
  solutions.
\newblock {\em Annals of Applied Mathematics}, 39(3):239--263, 2023.

\bibitem{jiao2023gas}
Yuling Jiao, Di~Li, Xiliang Lu, Jerry~Zhijian Yang, and Cheng Yuan.
\newblock Gas: A gaussian mixture distribution-based adaptive sampling method
  for pinns.
\newblock {\em arXiv preprint arXiv:2303.15849}, 2023.

\bibitem{jiao2022rate}
Yuling Jiao, Yanming Lai, Dingwei Li, Xiliang Lu, Fengru Wang, Jerry~Zhijian
  Yang, et~al.
\newblock A rate of convergence of physics informed neural networks for the
  linear second order elliptic pdes.
\newblock {\em Communications in Computational Physics}, 31(4):1272--1295,
  2022.

\bibitem{van1996weak}
Aad~W Van Der~Vaart, Adrianus~Willem van~der Vaart, Aad van~der Vaart, and Jon
  Wellner.
\newblock {\em Weak convergence and empirical processes: with applications to
  statistics}.
\newblock Springer Science \& Business Media, 1996.

\bibitem{Anthony1999Neural}
Martin Anthony and Peter~L. Bartlett.
\newblock Neural network learning: Theoretical foundations.
\newblock {\em Ai Magazine}, 22(2):99--100, 1999.

\end{thebibliography}
    \bibliographystyle{unsrt}

\end{document}